\documentclass[11pt]{article}
 
\input psfig.sty
\psfull

\begin{document}
 
\newtheorem{lemma}{Lemma}
\newtheorem{theorem}{Theorem}
\newtheorem{corollary}{Corollary}
\newtheorem{definition}{Definition}
\newtheorem{example}{Example}
\newtheorem{proposition}{Proposition}
\newtheorem{condition}{Condition}

\newcommand{\be}{\begin{equation}}
\newcommand{\ee}{\end{equation}}
\newcommand{\bea}{\begin{eqnarray}}
\newcommand{\eea}{\end{eqnarray}}
\newcommand{\beaa}{\begin{eqnarray*}}
\newcommand{\eeaa}{\end{eqnarray*}}
\newcommand{\ben}{\begin{enumerate}}
\newcommand{\een}{\end{enumerate}}
\newcommand{\bi}{\begin{itemize}}
\newcommand{\ei}{\end{itemize}}

\newcommand{\lip}{\langle}
\newcommand{\lan}{\langle}
\newcommand{\rip}{\rangle}
\newcommand{\ran}{\rangle}
\newcommand{\uu}{\underline}
\newcommand{\oo}{\overline}
\newcommand{\til}{\tilde}

\newcommand{\La}{\Lambda}
\newcommand{\la}{\lambda}
\newcommand{\eps}{\epsilon}
\newcommand{\vph}{\varphi}
\newcommand{\al}{\alpha}
\newcommand{\bet}{\beta}
\newcommand{\gam}{\gamma}
\newcommand{\Gam}{\Gamma}
\newcommand{\kap}{\kappa}
\newcommand{\Del}{\Delta}
\newcommand{\Th}{\Theta}
\newcommand{\s}{\sigma}
\newcommand{\sig}{\sigma}
\newcommand{\Sig}{\Sigma}
\newcommand{\del}{\delta}
\newcommand{\om}{\omega}
\newcommand{\Om}{\Omega}
\newcommand{\X}{{\bf \Xi}}

\newcommand{\integers}{Z\!\!\!Z}
\newcommand{\Z}{Z\!\!\!Z}
\newcommand{\N}{{\rm I\!N}}
\newcommand{\rationals}{{\rm I\!Q}}
\newcommand{\reals}{{\rm I\!R}}
\newcommand{\R}{{\rm I\!R}}
\newcommand{\RJ}{{{\rm I\!R}^J}}
\newcommand{\RJP}{{{\rm I\!R}^J_+}}
\newcommand{\naturals}{{\rm I\!N}}

\newcommand{\calA}{{\cal A}}
\newcommand{\calB}{{\cal B}}
\newcommand{\calC}{{\cal C}}
\newcommand{\calD}{{\cal D}}
\newcommand{\calF}{{\cal F}}
\newcommand{\calG}{{\cal G}}
\newcommand{\calH}{{\cal H}}
\newcommand{\calJ}{{\cal J}}
\newcommand{\calL}{{\cal L}}
\newcommand{\calM}{{\cal M}}
\newcommand{\calP}{{\cal P}}
\newcommand{\calS}{{\cal S}}
\newcommand{\calT}{{\cal T}}
\newcommand{\calU}{{\cal U}}
\newcommand{\calV}{{\cal V}}
\newcommand{\calX}{{\cal X}}
\newcommand{\calY}{{\cal Y}}

\newcommand{\proof}{\noindent {\bf Proof:\ }}
\newcommand{\proofOf}[1]{\noindent {\bf Proof of #1:\ }}
\newcommand{\remark}{\noindent {\bf Remark:\ }}
\newcommand{\remarks}{\noindent {\bf Remarks:\ }}
\newcommand{\note}{\noindent {\bf Note:\ }}
\newcommand{\esssup}{{\rm ess}\sup}
\newcommand{\essinf}{{\rm ess}\inf}
\newcommand{\pl}{\partial}
\newcommand{\noi}{\noindent}
\newcommand{\goto}{\to}
\newcommand{\ink}{\rule{.5\baselineskip}{.55\baselineskip}}
\newcommand{\qed}{\rule{.5\baselineskip}{.55\baselineskip}}

\def\ve{\varepsilon}
\def\vr{\varrho}

\newcommand{\diag}{{\rm diag}}
\newcommand{\trace}{{\rm trace}}
\newcommand{\tr}{{\rm tr}}
\newcommand{\w}{\wedge}
\newcommand{\dint}{\int\!\!\!\int}
\newcommand{\lt}{\left}
\newcommand{\rt}{\right}
\newcommand{\dist}{{\rm dist}}

\newcommand{\policy}{{u}}
\def\OBdry{\partial_o}
\def\CBdry{\partial_c}
\newcommand{\Sfrac}[2]{{{#1}\slash {#2}}}

\title{
An escape time criterion for queueing networks:
Asymptotic risk-sensitive control via differential games%
\thanks{This
research  was supported in part by the United States---Israel Binational
Science Foundation (BSF 1999179)}}
\author{Rami Atar\footnote{Electrical Engineering,
Technion---Israel Institute of Technology,
Haifa 32000, Israel. Research of this author also supported
in part by the fund for promotion of research at the Technion.},
Paul Dupuis\footnote{Lefschetz Center for Dynamical Systems,
Brown University,
Division of Applied Mathematics,
Providence,  R.I.\  02912. 
Research of this author also supported in part by
the National
Science Foundation (NSF-DMS-0072004, NSF-ECS-9979250) and the Army Research Office 
(DAAD19-99-1-0223).}
\ and Adam Shwartz\footnote{Electrical Engineering,
Technion---Israel Institute of Technology,
Haifa 32000, Israel. Research of this author also supported
in part by INTAS grant 265, and in part by the fund for promotion
of research at the Technion.}\\[.2in]
 }

\date{February 2, 2003}

\maketitle
 
\begin{abstract}
We consider the problem of risk-sensitive control of a stochastic network.
In controlling such a network,
an escape time criterion can be useful if one wishes
to regulate the occurrence of large buffers and buffer overflow. 
In this paper a risk-sensitive escape time criterion is formulated,
which in comparison to the ordinary escape time criteria penalizes exits which occur on short time intervals more heavily.
The properties of the risk-sensitive problem are studied in the large buffer limit,
and related to the value of a deterministic differential game with constrained dynamics.
We prove that the game has value, and that the value
is the (viscosity) solution of a PDE.
For a simple network, the value is computed,
demonstrating the applicability of the approach.
\end{abstract}

\section{Introduction}

In this paper we consider a problem of risk-sensitive control
(or rare event control) for queueing networks.  
The network includes
servers that can offer service to two or more
classes of customers, and a choice must be made
regarding which classes to offer service at each time.
We study a stochastic control
problem in which this choice is regarded as the control, and where
the cost is a risk-sensitive version of the time to
escape a bounded set. Hence, fixing $c>0$,
and denoting by $\s$ the time when the queue-lengths process
first exits a given domain, we consider $E_xe^{-c\s}$ as a criterion
to be minimized.
Such a criterion penalizes short exit times more heavily
than ordinary escape time criteria
(such as $E_x\s$, a criterion to be maximized).
There are at least two motivations for the use of such criteria when designing policies for the control of a network.
The first is that in many communication networks system performance is measured in terms of rare event probabilities (e.g., probabilities of data loss or excessive delay).
The second motivation follows from the connection between risk-sensitive controls and robust controls.
Indeed, as discussed in \cite{dupjampet},
the optimization of a single fixed stochastic network with respect to a risk-sensitive cost criteria automatically produces controls with specific and predictable robust properties.
In particular,
these controls give good performance for a family of perturbed network models (where the perturbation is around the design model and the size of the perturbation is measured by relative entropy), and with respect to a corresponding ordinary (i.e., not risk-sensitive) cost.

In many problems,
one considers the limit of the risk-sensitive problem as a scaling parameter of the system converges,
in the hope that the limit model is more tractable.
We follow the same approach here,
and show that the normalized costs in the risk-sensitive problems converge to the value function of a differential game with constraints.
As is well known, 
the convergence analysis is closely related to the large deviation properties of the sequence of controlled processes.
An interesting feature in the setting of stochastic networks is that the asymptotic analysis of a sequence of {\it controlled} networks is in many ways simpler than the analogous asymptotic analysis of a sequence of {\it uncontrolled} networks.
For example,
if one were to fix a particular state feedback service policy at each station,
then the calculation of the large deviation asymptotics is very difficult.
In contrast,
it turns out that calculation of the large deviation asymptotics of the optimally controlled network is quite feasible.
This is largely due to the fact that a {\it fixed} service policy invariably includes some {\it state discontinuities}.
For example,
a priority policy switches drastically when the highest priority queue empties.
When the policy is left as a parameter that is to be optimized these sharp discontinuities are not dealt with directly,
since the control and the large deviation behavior are identified simultaneously.
The situation is analogous to one found in the control of unconstrained processes such as diffusions.
If a fixed nonsmooth feedback control is considered then large deviation asymptotics are generally intractable,
but when the combined large deviation and optimal control problem is considered, much is possible \cite{dupmce}.

For simplicity, we restrict in this paper to a class of Markovian networks,
and consider just one simple cost structure.
Much more general statistical models can be treated with similar arguments,
as can a more general cost.  
A more fundamental restriction is on the routing in the network.
We assume a re-entrant line structure, so 
that the input streams follow a fixed route through the network--we do not allow either randomized or controlled routing.
Relaxing the last conditions
leads to a problem that is significantly more difficult to analyze,
and would require a considerable extension of the results we prove.

The deterministic game that is associated with the limit stochastic control
problem involves two players. One player
allocates service in a way analogous to the control in
the stochastic control problem, and the other player perturbs
the service and arrival rates. The cost is expressed in terms
of the large deviation rate function for the underlying arrival and service 
processes, cumulated up to the time the dynamics exit the domain.
Heuristically, the first player
identifies those classes it is most worthwhile to allocate
service to, so as to delay the escape as much as possible and thereby maximize the cost.
The player who selects the perturbed rates attempts to
minimize the cost by driving it out of the domain,
while paying a cost for perturbing the rates.

Our main result states that as the scaling parameter of the system
converges, the value for the stochastic control problem
converges to the value of the game. By way of proving the
result, we also show that the Hamilton-Jacobi-Bellman
equation associated with the game has a unique (Lipschitz continuous)
viscosity solution.

Several works have considered
problems of optimal exit probabilities in the context of controlled
diffusion processes, in the asymptotically small white noise
intensity regime.
Fleming and Souganidis \cite{flesou} use viscosity
solutions techniques to study a controlled diffusion where the
control enters in the drift coefficient. Dupuis and Kushner \cite{dupkus}
extend their results to
the case where the diffusion coefficient is possibly degenerate.
Their technique relates the stochastic control problem
to the game in a more direct way, using time discretization,
without involving PDE analysis.
The stochastic control problem studied in the current paper
has the property that the jump rates in certain directions
(those that correspond to services, not to arrivals)
can be controlled to assume arbitrarily low values, including zero.
It appears to be a more subtle problem than the ones in the above
cited papers, in that it is analogous to a
controlled diffusion problem where the control
enters also in the {\em diffusion} coefficient, and where no uniform
non-degeneracy condition is assumed. This kind of degeneracy makes
it difficult to apply the time discretization idea of \cite{dupkus}.
The main idea of \cite{dupkus}, in which one directly relates
the control problem to the game, is still fruitful in the
current setting. Following this approach,
we relate the limit inferior [resp., superior] of the 
asymptotic value for the control problem to the upper [resp., lower]
value of the game. However, showing that the game has value
and thereby obtaining the full convergence
result for the control problem requires a PDE analysis.

The PDE analysis uses viscosity solutions methods.
There are three types of boundary conditions associated with the PDE:
Neumann, Dirichlet, and ``state space constraint.''
The first two types of boundary conditions correspond in the game
to the nonnegativity constraint on queue lengths and to stopping upon exit from the domain, respectively.
The third type of boundary condition arises when there are portions
of the boundary where exit can be blocked unilaterally by
one of the players, and it is optimal for it to do so.
It is well known since Soner \cite{son} that such a scenario 
leads to the last boundary condition mentioned above.
Combining techniques of \cite{atadup2} and \cite{son},
we prove uniqueness of viscosity solutions for the PDE and show that
the game's upper and lower values are viscosity solutions,
thus establishing existence of value.
The trivial but crucial fact used in the uniqueness proof is
that the Isaacs condition holds (equation (\ref{eq:isaacs})).

As an example, we analyze a case
where the domain is a hyper-rectangle,
and where the network consists of one server and many
queues, each customer requiring service only once.
We find an explicit solution to the corresponding PDE,
assuming the parameter $c$ is large enough.
This is only an initial result in this direction,
but it shows that explicit solutions can be found.
The solution turns out to be of particularly simple form
(see equation (\ref{eq:form})). The optimal service discipline
stemming from the solution corresponds to giving priority
to class $i$ whenever the state of the system is within
a subset $G_i$ of the domain. The partitioning of the domain
into subsets has a simple structure too (see Figure \ref{fig:example}
in Section \ref{sec:example} for an example in two dimensions).
See \cite{atadupshw} for explicit solutions in the case of tandem queues,
as well as identities relating the perturbed rates with
the unperturbed ones in a more general network.

There is relatively little work on risk-sensitive and robust control of networks.
Ball et.\ al. have considered a robust formulation for network problems arising in vehicular traffic \cite{balday1,balday2},
and have explicitly identified the value function in certain instances.
Although their model is similar to ours in that the network dynamics are modeled via a Skorokhod Problem,
many other features, most notably the cost structure, are qualitatively different.
In addition,
the model they consider is not naturally related to a risk-sensitive control problem for a jump Markov model of a network.

The organization of the paper is as follows.
Section \ref{sec:setting} introduces the network and the stochastic
control problem, describes a key tool in our analysis, namely
the Skorokhod Problem (SP), introduces the
differential game, and states the main result.
Section~\ref{sec:limit} establishes
the relation between the control problem asymptotics
and the game's upper and lower values.
Characterization of the upper and lower values of the game
as viscosity solutions of a PDE, as well as uniqueness for this PDE
are established in Section \ref{sec:pde}.
Section~\ref{sec:example} presents an example, and the paper concludes with Section~\ref{sec:lemmas}, which gives the proofs of several lemmas.
Throughout the paper,
numbering such as Lemma $a.b$ refers to the $b$th item of Lemma $a$.

\section{Problem setting and the main result}\label{sec:setting}

\noi\uu{\bf The queueing network control problem.}
We consider a system with $J$ customer classes, and without loss 
assume that each class is identified with a queue at one of $K$ servers.
Each server provides service to at least one class.
Thus if $C(k)$ denotes the set of classes that are served by $k$,
then the control determines who receives service effort at server $k$
from among $i \in C(k)$. In particular, the sets $C(k)$, $k=1,\ldots, K$ are disjoint, 
with $\cup_k C(k)=\{1,\ldots,J\}$.
The state of the network is the vector of queue lengths, denoted
by $X$.
After a customer of class $i$ is served, it turns into a customer of
class $r(i)$, where $i=0$ is used to denote the ``outside.''
We let $ e_j $ denote the unit vector in direction $j$ and set $ e_0 = 0 $
so that following service to class $j$ the state changes by
$e_{r(j)} - e_j$.
The control will be described by the vector $u=(u_1,...,u_J)$,
where $u_i=1$ if class $i$ customers  are given service and $u_i=0$
otherwise. Since service can be given at any moment
to only one class at each station, the control vector must satisfy $\sum_{i\in C(k)}u_i\le1$ for
each $k$.
We next consider the scaled process $X^n$ under the scaling
which accelerates time by a factor of $n$ and shrinks space by
the same factor.
We are interested in a
risk-sensitive cost functional that is associated with
exit from a bounded set.
Let $G$ be a bounded subset of $\RJP $ that contains the origin
(additional assumptions on $G$ are given in Condition \ref{cond:G}).
Define
\[
\sigma^n \doteq \inf \{t:X^n(t) \not \in G \}.
\]
Then the control problem is to minimize the cost
$E_xe^{-nc\sigma^n}$, 
where $E_x$ denotes expectation starting
from $x$, and $c>0$ is a constant.
With this cost structure ``risk-sensitivity'' means that atypically short exit times are weighted heavily by the cost.
A ``good'' control will avoid such an event with high
probability. 
The significance from the point of view of stabilization of the
system is clear.
(See also \cite{dupmce} for the robust interpretation).

A precise description of the stochastic control problem is as follows.
Let $G^n \doteq n^{-1}\Z_+^J \cap G$. 
Define
$$
 U \doteq \left\{(u_i), i=1,\ldots,J:\sum_{i\in C(k)} u_i\le1, k=1,\ldots,K,
u_i\ge0, i=1,\ldots,J\right\}.
$$
For $u\in U$ and $f:\Z^J\to\R$ let
\be\label{eq:gen_constr}
\til{\cal L}^{u}f(x)
 \doteq \sum_{j=1}^{J}\lambda _{j}\left[ f(x+e_{j})-f(x)\right]
+\sum_{j=1}^J u_j \mu_j 1_{\{x+\tilde v_j\in\Z_+^J\}}
 \lt[ f(x+\tilde v_j)-f(x)\right],
\ee
where $\til v_j=e_{r(j)}-e_j$.
It is assumed that for each $i$, $\la_i\ge0$,
while $\mu_i>0$.
For each $n\in\N$ consider the scaling defined by
\be\label{eq:gens}
\til\calL^{n,u}f(x) \doteq n\til\calL^ug(nx),
\ee
where $f:n^{-1}\Z^J\to\R$ and $g(\cdot)=f(n^{-1}\cdot)$.
A {\em controlled Markov process} starting from $x$ will consist of a 
complete filtered
probability space $(\Om,\calF,(\calF_t),P_x^{n,u})$,
a state process $X^n$ taking values in  $G^n$ that is continuous from the right
and with limits from the left,
a control process $u$ taking values in  $U$,
such that $X^n$ is adapted to $\calF_t$,
$u$ is measurable and adapted to $\calF_t$,
$P_x^{u,n}(X^n(0)=x)=1$,
and for every function $f:n^{-1}\Z^J\to\R$
\[
f(X^n(t))-\int_0^t\til\calL^{n,u(s)}f(X^n(s))ds
\]
is an $\calF_t-$martingale.
$E_x^{n,u}$ denotes expectation with respect to $P_x^{n,u}$.
For a parameter $c>0$, the value function for the stochastic control problem is defined by
\be\label{eq:control}
  V^n(x) \doteq -\inf n^{-1}\log E_x^{u,n}e^{-nc\sig_n},
  \quad x\in G^n.
\ee
In this definition the infimum is over all controlled Markov processes.

A measurable function $u(x,t)$, $u:G^n\times[0,\infty)\to U$
is said to be a {\em feedback control}.
We will make use of two well known facts:
to each feedback control there corresponds a controlled Markov process with $u(t)=u(X^n(t),t)$,
and in the definition of the value function the infimum can be restricted to feedback controls.

In the formulation just given we allow the maximizing player to choose a control from the convex set $U$.
This is a relaxed formulation,
which allows the server to simultaneously split the effort between 2 or more customer classes.
An alternative control space that is more natural in implementation consists of only the vertices of $U$,
in which case the server can only server one class at a time.
In a general game setting,
the distinction between such ``relaxed'' and ``pure'' control spaces can be significant.
However,
in the present setting it will turn out that the value is the same for both cases.
This is essentially due to the fact that the game arises from a risk-sensitive control problem,
which imposes additional structure on the game,
and will be further commented on below.

\noi\uu{\bf Dynamics via the Skorokhod Problem.}
Our main goal will be to study the asymptotics of $V^n$,
and in particular, to show that they are
governed by the value of a deterministic differential
game. In order to define the dynamics of this game we first need a
formulation of the {\em Skorokhod Problem} (SP). We 
give here the simplest formulation which covers our needs.
The reader is referred to \cite{dupram23} for a more general framework.
Let 
\[
 D_{+}([0,\infty ):\RJ ) \doteq \left \{\psi\in D([0,\infty ):\RJ ):
   \psi(0)\in \RJP \right \},
\]
where $D([0,\infty ):\RJ )$ is the space of left continuous functions with right
hand limits, endowed with the uniform on compacts topology.
When restricting to continuous functions
we replace ``$D$'' with ``$C $''.
Let a set of vectors $\{ \gamma_i , \ i=1, \ldots , J \}$
be given and set $I(x) \doteq \{i:x_i=0\}$. 
For each point $x$ on $ \partial \RJP $ -- the boundary of $\RJP $ -- let
\[
d(x) \doteq \left\{ \sum_{i \in I(x)}a_i \gamma_i: a_i \geq 0, \left\|
\sum_{i \in I(x)}a_i \gamma_i \right\| = 1 \right\}.
\]
The Skorokhod Map (SM) assigns to every path $\psi\in
D_{+ }([0,\infty ) :\RJ )$ a path $\phi$ that starts at $\phi(0)=\psi(0)$,
but is constrained to $\RJP $ as follows.
If $\phi$ is in the interior of $ \RJP $ then the evolution of $\phi$
mimics that of $\psi$,
in that the increments of the two functions are the same until $\phi$ hits
$ \partial \RJP $.
When $\phi$ is on the boundary a constraining ``force'' is applied to keep
$\phi$ in the domain,
and this force can only be applied in one of the directions $d(\phi(t))$,
and only for $t$ such that $\phi(t)$ is on the boundary.
The precise definition is as follows.  
For $\eta\in
D([0,\infty ) : \RJ )$ and $t\in  [0, \infty)$ we let $|\eta|(t)$ denote the total
variation of
$\eta$ on $[0,t]$ with respect to the Euclidean norm on $ \RJ $.

\begin{definition}\label{def:sp}
Let $\psi\in D_{+}([0,\infty ) : \RJ )$ be given. Then $(\phi,\eta)$ solves
the SP for $\psi$ (with respect to $ \RJP $ and $\gamma_i, i=1,...,J $) if
$\phi(0)=\psi(0)$, and
if for all $t\in[0,\infty ) $
\begin{enumerate}
\item $\phi(t)=\psi(t)+\eta(t)$,
\item $\phi(t)\in \RJP $,
\item $|\eta|(t)<\infty$,
\item $|\eta|(t)=\int_{[0,t]} 1_{\{\phi(s)\in\partial  \RJP \}}
d|\eta|(s)$,
\item There exists a Borel measurable function
      $\gamma : [0,\infty ) \to \RJP $ such that $d|\eta|$-almost everywhere
      $\gamma(t)\in d(\phi(t))$, and such that
$$ \eta(t) = \int_{[0,t]}\gamma(s)d|\eta|(s).
$$
\end{enumerate}
\end{definition}
Under a certain condition on $\{ \gamma_i \}$ (known in the literature as the
{\it completely}-${\cal S}$ condition \cite{reiwil}), it is known
that solutions to the SP exist in all of $D_{ +}([0,\infty): \RJP )$.
Under further conditions 
(namely, {\em existence of the set $B$} -- see~\cite{harrei,dupish1,dupram23} and also Lemma \ref{lem:mu} below),
it is known that the Skorokhod Map is Lipschitz continuous, and consequently
the solution is unique.
Denoting the map $\psi\mapsto\phi$ in Definition \ref{def:sp} by
$\Gamma$, the Lipschitz property states that there is a constant
$K_1 $ such that
\be\label{sp:lip}
 \sup_{t\in[0,\infty)}\|\Gamma(\psi_1)(t)-\Gamma(\psi_2)(t)\|
 \le K_1 \sup_{t\in[0,\infty)}\|\psi_1(t)-\psi_2(t)\|,\quad
 \psi_1,\psi_2\in D_{ + }([0,\infty): \RJ ).
\ee
The SP that will be considered here is the one for which
$\gamma_i=e_i-e_{r(i)} = - {\tilde v}_i $. For this problem, the following is well known.
\begin{theorem}[\cite{harrei,dupram23}]\label{th:SP}
The SP associated with the domain $ \RJP $ and the constraint vectors
$\gamma_i, i=1,...,J $
possesses a unique solution, and the Skorokhod Map is
Lipschitz continuous on the space $D_{ +}([0,\infty): \RJP )$.
Moreover, the Skorokhod Map takes 
$ C_{+}([0,\infty): \RJ ) $ into $ C_{+}([0,\infty): \RJ ) $,
and therefore $ \Gamma (\phi ) $ is continuous if $ \phi $ is.
\end{theorem}
We next define a constrained ordinary differential equation.
As is proved (in greater generality) in \cite{dupish1},
one can define a projection $\pi: \RJ  \goto  \RJP $ that is consistent with
the constraint directions $\{\gamma_i, i=1,...,J \}$,
in that $\pi(x) = x$ if $x \in  \RJP $,
and if $x \not \in  \RJP $ then $\pi(x)-x = \alpha r$,
where $\alpha \geq 0$ and $r \in d(\pi(x))$.
With this projection given,
we can now define 
for each point $x \in \partial  \RJP $ and each $v \in  \RJ $
the {\it projected velocity}
\[
\pi(x,v) \doteq \lim_{\Delta \downarrow 0} \frac{\pi(x+\Delta
v)-\pi(x)}{\Delta}.
\]
For details on why this limit is always well defined and further properties
of the projected velocity we refer to \cite[Section 3 and Lemma
3.8]{buddup}.
Let $v:[0,\infty) \goto  \RJ $ have the property  that each component of
$v$ is integrable over each interval $[0,T]$, $T<\infty$.
Then the ODE of interest takes the form
\be
\label{eqn:ODE}
\dot \phi(t) = \pi(\phi(t),v(t)), \quad \phi (0) = \phi_0 \in \RJP .
\ee
An absolutely continuous function $\phi:[0,\infty) \goto  \RJP $ is a solution
to (\ref{eqn:ODE}) if  the equation is satisfied a.e.\ in $t$.
By using the regularity properties (\ref{sp:lip}) of the associated Skorokhod Map
 and because of the
particularly simple nature of the right hand side,
one can show that $\phi$ solves (\ref{eqn:ODE}) if and only if $\phi$ is
the image of $\psi(t) \doteq \int_0^t v(s)ds +x$ under the Skorokhod
Map, and thus all the standard qualitative properties (existence and
uniqueness of solutions, stability with respect to perturbations, etc.)
hold \cite{dupish1,dupnag}.

As mentioned above, the SP formulation will be our means of defining
the dynamics of a deterministic game. Before discussing this game,
let us show how the same formulation is also useful for
the stochastic control problem defined earlier in this section.
First, since $\til v_j = e_{r(j)}-e_j = - \gamma_j$, it is easy to verify that for the particular SP considered here
$\pi(x,v)=v1_{x+v\in\Z_+^J}$ for all $x\in\Z_+^J$ and
$v\in\{\til v_j:j=1,\ldots,J\}$.
Therefore the generator $\til\calL^{u}$ of (\ref{eq:gen_constr})
can also be written as
$$
 \til\calL^uf(x)=\sum_{j=1}^J\la_j[f(x+e_j)-f(x)]+
 \sum_{j=1}^Ju_j\mu_j[f(x+\pi(x,\tilde v_j))-f(x)].
$$
A measurable function $u(t)$, $u:[0,\infty)\to U$
is said to be an {\em open loop control}. 
Note that this control has no state feedback.
When $u$ is an open loop control,
it is possible to write the corresponding controlled process $X$ as $\Gamma(Y)$.
The process $Y$, which will be called
{\em the unconstrained controlled process},
is a controlled Markov process with a simpler structure.
To be precise, let
$$
 \til{\cal L}_0^{u}f(x)
 =\sum_{j=1}^{J}\lambda _{j}\left[ f(x+e_{j})-f(x)\right]
 +\sum_{j=1}^Ju_j\mu_j\left[ f(x+ \tilde v_j)-f(x)\right],
$$
and let $\til\calL_0^{n,u}$ be defined analogously to (\ref{eq:gens}).
A controlled Markov process
$X^n$ on $G^n$ [respectively, $Y^n$ on $n^{-1}\Z^J$] is defined as before, 
but now using the
generator $(\calL^{n,u}f)(t)=(\til\calL^{n,u(t)}f)(x)$
[resp., $(\calL_0^{n,u}f)(t)=(\til\calL_0^{n,u(t)}f)(x)$].
The simplification that the SP introduces
is that if $u$ is an open loop control, and if $Y^n$ is a
controlled Markov process corresponding to $\calL_0^{n,u}$
on $(\Om,\calF,(\calF_t),P_x^{n,u})$, then $X^n=\Gamma(Y^n)$
is a controlled Markov process corresponding to $\calL^{n,u}$
on the same filtered probability space.
The role played by the SP in relating constrained and
unconstrained processes is exhibited here in a simple fashion,
for introductory purposes. We will, in fact, use it in
a slightly more complicated setting later on in Lemma \ref{lem:xbar} and Lemma 
\ref{lem:xbar2}.

\noi\uu{\bf A differential game.}
In this paper,
we prove that the value function $V^n(x)$ for a stochastic control problem associated with our queueing network model is approximately equal (for large $n$) to the value function of a related differential
game.
In addition,
the dynamics of this game are defined in terms of an associated SP. 
Before introducing the game
formally we explain why this is to be expected.
In a problem with no control, the exponential decay rate
of quantities such as $Ee^{-cn\sig_n}$ is
given in terms of the sample path
large deviation rate function associated with the process, which in
turn can be expressed in terms of the rate function for the Poisson
primitives that drive the model.
This is supported by the
well known Laplace's principle \cite{dupell}. Heuristically,
one thinks of the rate function as a cost paid for changing the measure
so as to make the rare event of exiting
on short time interval a probable event.
Laplace's principle asserts that the decay rate can be expressed
as the solution to a deterministic optimization problem
involving the cost $-c\sig$ combined with the cost of changing the measure:
cf.\ \cite[Eq.\ (5.20)--(5.23)]{shwe}.
When the stochastic model involves optimal control, there is one
more variable to optimize over in the limit, and this results in a game.
The game's deterministic dynamics are the natural law of large numbers limit under the changed measure.
Boundary constraints and constraining meachanisms which are present in the prelimit model are represented in the limit model by the SP.
The cost for the game involves
the large deviation rate function for the Poisson primitives,
the time till the dynamics exits the domain, and the parameter $c$.

We thus consider a zero sum game involving two players.
One (which we call the maximizing player)
selects the service allocation and attempts to maximize.
The other (called the minimizing player)
chooses the perturbed arrival and service rates and attempts to minimize.
Throughout, the perturbed rates will be denoted by an overbar,
as in $\bar\la_i,\bar\mu_i$.

Recall that for $u\in U$, $u_i$
stands for the fraction of service effort given to class $i$.
The control space for the maximizing player is
$$ \bar U\doteq \{u:[0,\infty)\to U\ ;\ u \mbox{ is measurable}\}.
$$
Let $l:\R\to\R_+\cup\{+\infty\}$
be defined as
$$
 l(x)\doteq \lt\{\begin{array}{ll}x\log x-x+1 & x\ge0,  \\ +\infty & x<0,
     \end{array}\rt.
$$
where $0\log0\doteq0$.
Denoting $M =[0,\infty)^{2J}$, the control space for the minimizing
player will be
\begin{equation}
\label{embar}
 \bar M =\{m=(\bar\la_1,\ldots,\bar\la_J,\bar\mu_1,\ldots,\bar\mu_J):
[0,\infty)\to M;
\ \mbox{$m$ is measurable, $l\circ m$ is  locally integrable}\}.
\end{equation}
For $u\in U$ and $m\in M$ define
\[
v(u,m) \doteq \sum_{j=1}^J \bar \lambda_jv_j + \sum_{i=1}^J
u_i \bar\mu_i \tilde v_i,
\]
where $v_j=e_j$, and as before $\tilde v_i=e_{r(i)}-e_i$.
Then the dynamics are given by
\[
\lt\{\begin{array}{ll} \dot \phi(t) = \pi (\phi(t),v(u(t),m(t))), \\ 
 \phi(0)=x.
 \end{array}\rt.
\]
To define the cost for the game, let
$\rho:U\times M\to\R_+\cup\{+\infty\}$ be
$$ \rho(u,m)\doteq \sum_{i=1}^{J}\lambda _{i}l\left( \frac{\bar{\lambda}_{i}}{
\lambda _{i}}\right) +\sum_{i=1}^{J}u_i\mu _il\left( \frac{\bar{\mu}_i}
{\mu _i}\right).
$$
By convention, if $\la_i=0$ and $\bar\la_i>0$ for some $i$, we let
$\rho=\infty$ (recall that by assumption, $\mu_i>0$).
Let the exit time be defined by
$$
 \sig\doteq \inf\{t:\phi(t)\not\in G\}.
$$
With $c > 0$ as in~(\ref{eq:control}),
the cost is given by
$$ C(x,u,m)=\int_0^\sig [c+\rho(u(t),m(t))]dt.
$$

As in \cite{ellkal} we need the notion
of strategies. We endow both $\bar U$ and $\bar M$ with the metric
$\til\rho(\om_1,\om_2)=\sum_n2^{-n}(\int_0^n
|\om_1(t)-\om_2(t)|dt \w1)$,
and with the corresponding Borel $\s$-fields.
A mapping $\al:\bar M \to\bar U$ is called a
{\em strategy for the maximizing player}
if it is measurable
and if for every $m,\tilde m\in\bar M$ and $t>0$ such that
$$ m(s)=\tilde m(s) \mbox{ for a.e.\ } s\in[0,t],
$$
one has
$$ \al[m](s) = \al[\tilde m](s) \mbox{ for a.e.\ } s\in[0,t].
$$
In an analogous way, one defines a mapping $\beta:\bar U\to \bar M$
to be a {\em strategy for the minimizing player}.
The set of all strategies for the maximizing [resp., minimizing]
player will be denoted by $A$ [resp., $B$].
The lower value for the game is defined as
\[
V^-(x) = \inf_{\beta\in B} \sup_{u\in\bar U}C(x,u,\beta[u]),
\]
and the upper value as
$$
V^+(x)=\sup_{\al\in A}\inf_{m\in\bar M}C(x,\al[m],m).
$$
To avoid confusion, we remark that despite the terms ``upper'' and
``lower'' value,  it is not in general obvious that $V^-\le V^+$.

\noi\uu{\bf Main result.}
We make the following assumption on the domain $G$. Let
$$
 \calJ_+\doteq\{i\in\{1,\ldots,J\}:\la_i>0\}.
$$

\begin{condition}\label{cond:G}
We assume that the domain $G$ satisfies one of the following.
\begin{enumerate}
\item
$G$ is a rectangle given by
$$ G=\{(x_1,\ldots,x_J): 0\le x_i<z_i, i\in\calJ_+;\ 0\le x_j\le z_j,
j\not\in\calJ_+\},
$$
for some $z_i>0$, $i=1,\ldots,J$.
\item
$G$ is simply connected and bounded, and given by
$$
 G=\bigcap_{i\in\calJ_+}G_i,
$$
where for $i\in\calJ_+$, we are given positive Lipschitz functions
$\phi_i:\R^{J-1}\to\R$, and
$$
 G_i=\{(x_1,\ldots,x_J)\in \RJP : 0\le x_i<\phi_i(x_1,\ldots,x_{i-1},
x_{i+1},\ldots,x_J)\}.
$$
\end{enumerate}
\end{condition}
This condition covers many typical constraints one would consider on buffer size, 
including separate constraints on individual queues (Condition \ref{cond:G}.1) and one constraint on the sum of the queues (Condition \ref{cond:G}.2).

\begin{figure} 
\centerline{
\begin{tabular}{cc}
\psfig{file=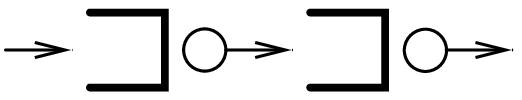}\qquad\qquad\qquad
\psfig{file=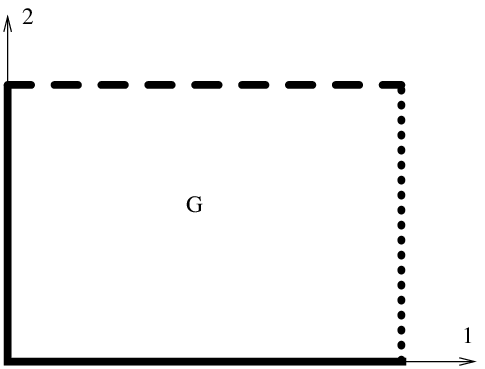}
\end{tabular}
}
\caption{
A simple queueing network, a rectangular domain and three
types of boundary.
Full line: $~$ \hfill $\pl_+G$, dashed line: $\pl_cG$,
and dotted line: $\pl_oG$
\hfill $~$ }
\label{fig1}
\end{figure}

The shape of the domain is simpler in   Condition \ref{cond:G}.1,
in that it is restricted to a hyper-rectangle.
On the other hand, it is also possible under this condition for the
maximizing player to unilaterally prevent an exit 
through a certain portion of $\pl G\setminus\pl\R_+^J$.
Although it is in principle possible that the dynamics
could exit through this portion of the boundary, it will always be optimal for
the maximizing player to not allow it.
Consider the simple network illustrated in Figure \ref{fig1}.
The maximizing player can prevent exit through the dashed portion of the boundary simply by stopping service at the first queue. 
As a consequence,
there are in general three different types of boundary--the constraining boundary due to non-negativity constraints on queue length,
the part of the boundary where exit can be blocked,
and the remainder.
These three types of boundary behavior will result,
in the PDE analysis, in three types of boundary conditions.
We now define the three portions of the boundary.
Under Condition \ref{cond:G}.1, let
$$
\pl_cG=\{(x_1,\ldots,x_J)\in G: x_j=z_j,\,\mbox{some}\,j\not\in\calJ_+\}.
$$
For notational convenience, we let $\pl_cG=\emptyset$ under
  Condition \ref{cond:G}.2. In both cases we then set
$$
\pl_oG=\pl G\setminus G, \qquad
\pl_+G=(G\cap\pl\R_+^J)\setminus\pl_cG.
$$
Note that in both cases, $\pl_cG$, $\pl_oG$ and $\pl_+G$
partition $\pl G$. Also, $\pl_cG\subset G$ while $\pl_oG\subset
G^c$. As usual, we will denote $G^o=G\setminus\pl G$
and $\bar G=G\cup\pl G$.
$ \pl_cG $ is the part of the boundary were the
maximizing player can prevent the dynamics from exiting,
and $ \pl_oG $ is the part where it can not.
Finally, it will be convenient to denote
$$
\pl_{co}G=\pl_oG\cup\pl_cG.
$$

Our main result is the following.
\begin{theorem}\label{th:main}
Let Condition \ref{cond:G} hold.
Then $V^+=V^-\doteq V$ on $G$. Moreover, if $x_n\in G^n$,
$n\in\N$ are such that $x_n\to x\in G$, then
$\lim_{n\to\infty}V^n(x_n) = V(x)$.
\end{theorem}

\remark 
A stronger form of the convergence statement in fact holds.
Namely,
$$
 \limsup_{\eps\downarrow0}
 \limsup_{n\to\infty}\sup\{|V^n(x)-V(y)|: x\in G^n,y\in G, |x-y|\le\eps\}
 =0.
$$
This is an immediate consequence of Theorem \ref{th:main}
and the fact that for each $n$, $V^n$ is Lipschitz on $G^n$,
with a constant that does not depend on $n$ (Lemma \ref{lem:unif}).
\qed

The proof is established in two major steps.
Step 1 will be an immediate consequence of the main results of
Section \ref{sec:limit},
and Step 2 will follow from Section \ref{sec:pde}.

\noi{\em Step 1.} We define a version of the game, technically
easier to work with, in which
all perturbed rates ($\bar\la_i,\bar\mu_i$) are bounded
by $b<\infty$. The corresponding upper and lower values,
defined analogously, are denoted by $V^{b,+}$ and $V^{b,-}$.
Then we show that for all $b$ large enough
(cf.\ Theorem \ref{th:limit})
$$
 V^{b,+} (x)\le\liminf_{n\to\infty}V^n(x_n)\le\limsup_{n\to\infty}V^n(x_n)
 \le V^{b,-} (x) .
$$

\noi{\em Step 2.}
We show that for $b$ large,
$V^{b,+} =V^{b,-} $ on $G$. To this end, we formulate a PDE
for which we show that uniqueness of (Lipschitz) viscosity solutions
holds (Theorem \ref{th:unique}),
and also show that both $V^{b,+}$ and $V^{b,-}$ are viscosity solutions
(Theorem \ref{th:solve}).
Since $ V^n (x) $ does not depend on $b$, neither do $ V^{b,\pm} (x) $.
Theorem \ref{th:main} follows.

\section{The control problem and the game}\label{sec:limit}

We begin by stating some basic properties of the stochastic
control problem and of the deterministic game.
The proofs of these properties are 
deferred to Section \ref{sec:lemmas}.

Consider the following generators, defined for any $u\in U$ and
$m\in M$, for constrained and unconstrained controlled Markov
processes:
\beaa
\calL^{n,u,m}f(x) &=& \sum_{j=1}^J n \bar \lambda_j \left [f \left (x+\frac{1}{n} v_j \right )-f(x) \right ]
+ \sum_{i=1}^J n \bar \mu_i u_i\left [ f \left ( x+\frac{1}{n}\pi(x, \tilde v_i) \right )-f(x) \right ], \\
\calL_0^{n,u,m}f(x) 
&=& \sum_{j=1}^J n \bar \lambda_j \left [ f \left (x+\frac{1}{n} v_j \right )-f(x) \right ]
+ \sum_{i=1}^J n \bar \mu_i u_i \left [f \left ( x+\frac{1}{n} \tilde v_i \right )-f(x) \right ].
\eeaa
The definition of the corresponding controlled processes will be made precise in Lemmas
\ref{lem:xbar} and \ref{lem:xbar2}.

Owing to the logarithmic transform in (\ref{eq:control}), one expects $V^n$ to satisfy an Isaacs equation \cite{flesou}.
In fact, $V^n$ satisfies 
\be\label{eq:dpe}
\lt\{\begin{array}{cc}
0=\sup_{u\in U} \inf_{m\in M} [\calL^{n,u,m}V^n(x) +c+\rho(u,m)],
  & x \in G^n \\ 
V^n(x)=0, & x \not \in G^n.
\end{array}\rt.
\ee
We comment that this is also the dynamic programming equation (DPE)
for an associated stochastic game
that is related to the deterministic game
via a law of large numbers scaling and limit,
and will not be further considered in this paper.

\begin{lemma}\label{lem:dpe}
The value function $V^n$ of (\ref{eq:control}) uniquely solves
the DPE (\ref{eq:dpe}).
\end{lemma}

The following lemma gives a key estimate on the value function.

\begin{lemma}\label{lem:unif}
Under Condition~\ref{cond:G},
$ V^n (x) $ satisfies the Lipschitz property
on $(n^{-1}\Z_+^J) \cap \bar G$
with a constant that does not depend on $n\in\N$.
Consequently, $ \sup_{n, x \in G^n} V^n (x) < \infty $.
\end{lemma}
We comment that the above result is, in general, not valid for $V^n$
on $n^{-1}\Z_+^J$, since $V^n$ changes abruptly near
the portion $\pl_cG$ of the boundary.

For each fixed $u\in U$, the mapping $m\rightarrow \rho(u,m)$,
when restricted to $\bar \mu_i$ such that $u_i>0$, is strictly convex with compact level sets.
We conclude that the infimum over $m$ in the DPE
is achieved, 
and denote such a point by $m^n(x,u)$. 
Although part 1 of the following lemma is not used elsewhere,
it indicates why the Isaacs condition should hold in (\ref{eq:dpe}).

\begin{lemma}\label{lem:bound:m} Let Condition \ref{cond:G} hold.
Then
\ben
\item
$m^n(x,u)$ can be chosen independently of $u$, and 
\item
there is $b_0<\infty$ such that for all $x,n$ and $u$, $ m^n(x,u) \leq b_0$. 
\een
\end{lemma}

We introduce two parametric variations of the game defined in
Section \ref{sec:setting}. The first will be associated with
domain perturbation (parameterized by the symbol $a$), and the second
with a bound on the perturbed rates (parameterized by the symbol $b$).

For some fixed $a_0>0$, consider perturbations $G_a$,
$a\in(-a_0,a_0)$ of the domain $G$
defined as follows. If $G$ satisfies Condition~\ref{cond:G}.1,
then $G_a$ is defined as $G$, but with $z_i$ replaced by $z_i+a$,
$i=1,\ldots,J$.
If $G$ is as in Condition~\ref{cond:G}.2, then $G_a$ is defined
as $G$, but where $\phi_i$ is replaced by $\phi_i+a$, $i\in\calJ_+$.

For any $b\in(0,\infty)$, let $M^b=[0,b]^{2J}$. Analogously to
the definition (\ref{embar}) of $\bar M$, let
\begin{equation}
\label{embarb}
 \bar M^b =\{m=(\bar\la_1,\ldots,\bar\la_J,\bar\mu_1,\ldots,\bar\mu_J):
[0,\infty)\to M^b \ ;\ \mbox{$m$ is measurable}\}.
\end{equation}
Strategies and values for the game are then defined
analogously to the way strategies and values are defined for
the original game, using $\bar M^b$ in place of $\bar M$.
It will be convenient to set $M^\infty\doteq M$ and
$\bar M^\infty\doteq\bar M$, and to refer to the original game
of Section \ref{sec:setting} as the case $b=\infty$.

The cost, sets of strategies, lower and upper values of the
games resulting by the introduction of the parameters
$a$ and $b$ will be denoted as $C_a(x,u,m)$, $A^b$, $B^b$,
$V^{b,-}_a$ and $V^{b,+}_a$. When $a=0$ [resp., $b=\infty$],
the dependence on $a$ [$b$] will be eliminated from the notation,
as in $V^-_a$, $V^{b,-}$.

Let $b_0$ be as in Lemma \ref{lem:bound:m}.
Denote
\begin{equation}
\label{eq:lal}
b^* \doteq \max \{ b_0, \lambda_i , \mu_i , i = 1 , \ldots , J \} + 1.
\end{equation}

\begin{lemma}\label{lem:cont}
Assume Condition~\ref{cond:G}.
Then
\ben
\item
$\dist(\pl_{co} G_a,\pl_{co} G)\doteq 
\inf\{|x-y|:x \in \pl_{co} G_a, y \in \pl_{co} G\} >0\mbox{ if } 0<|a|<a_0$;
\item
the values
$V^{b,\pm}$ are bounded on $G$, uniformly for $b\in[b^*,\infty]$,
and there is a constant $c_0$ such that for any
$x\in G$, $|a|<\eps$ (where $\eps$ depends on $x$), and $b\in[b^*,\infty]$,
one has
$|V^{b,-}_a(x)-V^{b,-}(x)|\le c_0|a|$ and
$|V^{b,+}_a(x)-V^{b,+}(x)|\le c_0|a|$.
\een
\end{lemma}

The following lemma shows that any nearly optimal strategy for the minimizing player will satisfy a uniform upper bound on the integrated running cost.
Moreover, there is a finite time $T_0$ such that 
for each such minimizing strategy, any open loop control used by the maximizing player leads to  exit by $T_0$.
Similarly,
given any strategy for the maximizing player the minimizing player can restrict to open loop controls that force exit by $T_0$.
\begin{lemma}\label{lem:apriori}
Fix $b\in[b^*,\infty]$.
Given $\beta \in B^b$,
write
$(\bar\la_i(\cdot ),\bar\mu_i(\cdot ))=\beta[u](\cdot)$.  For $z,T>0$ let
$B^{z,T}$ denote the set of $\beta\in B^b$ which satisfy 
$$
  \int_0^T\sum_i[\la_il(\bar\la_i(t)/\la_i)+u_i(t)
  \mu_il(\bar\mu_i(t)/\mu_i)]dt
  \le z,
$$
for all $u\in\bar U$.
For $\al\in A^b$, let $\bar M(\al,T)$ denote the set of $m\in\bar M$
for which $\s(x,\al[m],m)\le T$.
Then there are constants $z_0,T_0>0$ such that
$$
V^- (x)=\inf_{\beta\in B^{z_0,T_0}}\sup_{u\in\bar U}
\int_0^{\sig\w T_0}[c+\rho(u(t),\beta[u](t))]dt,
$$
and
$$
V^+(x)=\sup_{\al\in A}\inf_{m\in\bar M(\al,T_0)}
\int_0^{\sig\w T_0}[c+\rho(\al[m](t),m(t))]dt.
$$
\end{lemma}

In the rest of the section  the strategies $\beta$
will be assumed (without loss)
to be in $B^{z_0,T_0}$, where $z_0,T_0$ are as in
Lemma \ref{lem:apriori}, and are fixed throughout. Also, $m\in\bar M$
will be assumed to be in $\bar M(\al,T_0)$ whenever it is
clear which $\al$ is considered.
With an abuse of notation, we denote
$B^{z_0,T_0}$ by $B$.

\begin{lemma}\label{lem:vallip} Under Condition \ref{cond:G},
$V^{b,-}$ and $V^{b,+}$ are Lipschitz on $G$, uniformly for
$ b\in[b^*,\infty]$.
\end{lemma}

We are now ready to prove the following result.
\begin{theorem}\label{th:limit}
Let Condition \ref{cond:G} hold, and let $b\in[b^*,\infty)$.
Then for any $x\in G$ and any $\{ x_n \} $ converging to $x$ (with
$x_n\in G^n$),
$$
 V^{b,+}(x)\le\liminf_{n\to\infty}V^n(x_n)\le\limsup_{n\to\infty}
 V^n(x_n)\le V^{b,-}(x).
$$
\end{theorem}

\proofOf{Theorem~\ref{th:limit}}
The result is established by considering a sequence of stochastic
processes, defined using the constrained ODEs, but for which
the controls $u$ and $m$ are governed by, on one hand, a nearly optimal
strategy for the game, and on the other hand, a nearly optimal
control for the stochastic control problem.
The technique uses standard martingale estimates, and is based
on the construction (deferred to Section \ref{sec:lemmas})
of an auxiliary controlled Markov process
that is controlled by
the selected strategy and stochastic control.

\noi{\em \uu{Upper bound}}

\noi
Fix $b\in[b^*,\infty)$.
The dependence on $b$ will be suppressed in the notation for $V^-$, $V^-_a$, etc.
We first show that
\be\label{eq:ub}
\limsup_{n\to\infty} V^n(x_n) \leq V^-(x).
\ee
According to  Lemma~\ref{lem:cont}.2, it is enough to show that
for all $a>0$
\[
\limsup_{n\to\infty} V^n(x_n) \leq V_a^-(x).
\]
Let $\beta\in B^b $ and $a>0$ be given,
and set $C_a(x,\beta)=\sup_{u\in\bar U} C_a(x,u,\beta[u])$.  
It is enough to show that
\be\label{eq:ca}
\limsup_{n\to\infty} V^n(x_n) \leq C_a(x,\beta), \quad a>0.
\ee
We therefore fix $\beta$ throughout, and turn to prove (\ref{eq:ca}).
We can assume without loss that
\begin{equation}\label{eq:31}
C_a(x,\beta)<\infty.
\end{equation}
Note that in the DPE (\ref{eq:dpe}) the supremum is with respect to $u$ in
a compact set $U$, and that the function being maximized is continuous in $u$ (for each $y$).
Let $u^n(y)$ denote a point where it is achieved.
Then for any $m\in M$ and $y\in G^n$,
\be\label{ineq1}
0 \leq  \calL^{n,u^n(y),m}V^n(y) +c+\rho(u^n(y),m).
\ee
\begin{lemma}\label{lem:xbar}
Let $n$ be fixed, and let $b\in[b^*,\infty)$,
$\beta$ and $x_n$ be as above.
Then there is a filtered probability space
$(\bar\Om,\bar F,(\bar F_t),\bar P)$, and $\bar F_t$-adapted
RCLL processes
$\bar X^n$, $\bar Y^n$ and $m^n$ such that with $\bar P$-probability one
$m^n(t)=\beta[\bar u^n](t)$ a.e.\ $t$, $\bar u^n(t)=u^n(\bar X^n(t))$,
$\bar X^n=\Gamma(\bar Y^n)$,
$\bar X^n(0)=\bar Y^n(0)=x_n$,
and for any $f$
\[
f(\bar X^n(t)) - \int_0^t \calL^{n,\bar u^n, m^n(s)}f(\bar X^n(s))ds
\]
\[
f(\bar Y^n(t)) - \int_0^t \calL_0^{n,\bar u^n, m^n(s)}f(\bar Y^n(s))ds
\]
are $(F_t)$-martingales.
Moreover, with $T_0$ as in Lemma~\ref{lem:apriori} let $N_n $ denote
the total number of jumps of $ \bar Y^n $ on $ [0,T_0]$. Then
\be\label{ineq:Nn}
EN_n\le 2JT_0bn.
\ee
\end{lemma}

\proof See Section \ref{sec:lemmas}.

Returning to the proof of~(\ref{eq:ca}), let $\bar u^n(t) \doteq u^n(\bar X^n(t))$
and let $\bar\sigma^n$ be the first
exit time of $\bar X^n$ from $G$.
Combining
(\ref{ineq1}) and Lemma \ref{lem:xbar}, for any bounded stopping
time $S\le\bar\sig^n$,
\be\label{ineq2}
V^n(x_n) \leq \bar E_{x_n} \left[ V^n(\bar X^n(S))
+\int_0^{S}[c+\rho(\bar u^n(s),\beta[\bar u^n](s))]ds \right].
\ee
Denoting $\beta[\bar u^n](t)= \{ (\bar\la_i^n(t),\bar\mu_i^n(t)) \} $,
define $\phi^n$ as $\phi^n=\Gamma(\psi^n)$, where
\[
   \psi^n (t)= x+\int_0^t v(\bar u^n,\beta[\bar u^n])ds,
\]
and let
\[
\hat\sig^n_a \doteq \inf\{t:\phi^n (t) \not \in G_a\}.
\]
Then the definition of $C_a(x,\beta)$ implies
\[
\int_0^{\hat\sig^n_a}[c+\rho(\bar u^n(s),\beta[\bar u^n](s))]ds
 \leq C_a(x,\beta).
\]
Apply (\ref{ineq2}) with $S=\hat\sig_a^n\w\bar\sig^n\w T$.
If $T$ is sufficiently large, then (\ref{eq:31}) and the
fact that $c>0$ imply $\hat\sig_a^n \leq T$. Thus, using
$\bar E_{x_n}V^n(\bar X^n(\bar\sig^n)))=0$,
\[
V^n(x_n) \leq \bar E_{x_n} \left[ V^n(\bar X^n(\hat\sig^n_a))
1_{\{\hat\sig^n_a < \bar\sigma^n\}}
 +\int_0^{\hat\sig^n_a}[c+\rho(\bar u^n(s),\beta[\bar u^n](s))]ds \right].
\]
Again using the uniform boundedness of $V^n(x)$ (Lemma \ref{lem:unif}), there is a constant $b_2 <\infty$
such that for all $n$
\be\label{ineq3}
V^n(x_n) \leq b_2  \bar P_{x_n}(\hat\sig^n_a \le \bar\sigma^n)
 + C_a(x,\beta).
\ee
In what follows, we show that $\bar P_{x_n}(\hat\sig^n_a \le \bar\sigma^n)$
tends to zero. To this end,
note that $\calL_0^{n,u,m}\,{\rm id}(y)=\sum_i\bar\la_iv_i+\sum_iu_i\bar\mu_i
\tilde v_i$, where ${\rm id}$ is the identity map.
Therefore, using again Lemma \ref{lem:xbar},
\[
\bar Y^n(t) -x_n = \int_0^t \lt[ \sum_{i=1}^J \bar\la_i^n(s)v_i
 +\sum_{i=1}^J\bar u_i^n(s)\bar\mu_i^n(s)\tilde v_i\rt]
 ds + \eta^n(t),
\]
where $\eta^n$ is a zero mean martingale.
To prove that
\be\label{conv:mgale}
 \sup_{t\in[0,T]}|\eta^n(t)| \to 0 \quad \mbox{in distribution},
\ee
it is enough, by Doob's maximal inequality, to show that
$$
\bar E|\eta^n(T)|^2\to0.
$$
Let $[x](t)=\sum_{s\in[0,t]}|\Del x_s|^2$.
By the Burkholder-Davies-Gundy inequality (see \cite{delmey}, VII.92),
 $$
 \bar E|\eta^n(T)|^2 \le c_1\bar E[\eta^n](T),
$$
where $c_1$ is a constant.
Since each jump is bounded by $c_2n^{-1}$
($c_2$ a constant) and the total number of jumps
$N_n(T)$ satisfies (\ref{ineq:Nn}),
$$
\bar E|\eta^n(T)|^2 \le c_3n^{-2}\bar EN_n(T) \le c_4n^{-1},
$$
which proves (\ref{conv:mgale}).
This implies that $\sup_{[0,T]}|\bar Y^n(t)-\psi^n(t)|\to0$ in
distribution, and therefore the continuity of $\Gamma$ implies $\sup_{[0,T]}|\bar X^n(t)-\phi^n(t)|\to0$
in distribution.
By Lemma \ref{lem:cont}.1,
\beaa
\bar P_{x_n} \left( \hat\sig^n_a \le \bar\sigma^n \right)
&\leq &
\bar P_{x_n}(\bar X^n(\hat\sig_a^n)\in G,
\phi^n(\hat\sig_a^n)\in\pl_{co}G_a)\\
&\le&
\bar P_{x_n} \left (\sup_{t\in[0,T]}|\bar X^n(t)-\phi^n(t)|\ge b_1 \right ),
\eeaa
where $b_1>0$ depends only on $a$.
Hence by (\ref{conv:mgale}),
$\bar P_{x_n}[ \hat\sig^n_a \le \bar\sigma^n]\to0$
as $n\to\infty$. Therefore (\ref{ineq3}) implies
\[
\limsup_{n \goto \infty} V^n(x_n) \leq C_a(x,\beta).
\]
This gives (\ref{eq:ca}) and completes the proof of (\ref{eq:ub}).

\noi{\em \uu{Lower bound}}

\noi
Next we prove
\be\label{eq:lb}
\liminf_{n\to\infty} V^n(x_n) \geq V^+(x).
\ee
By Lemma \ref{lem:cont}.2, it is enough to show that
for all $a<0$
\[
\liminf_{n\to\infty} V^n(x_n) \geq V_a^+(x).
\]
Let $\al\in A$ be given,
and set $C_a(x,\al)=\inf_{m\in\bar M^b} C_a(x,\al[m],m)$.  
Then it suffices to show
\be\label{eq:ca2}
\liminf_{n\to\infty} V^n(x_n) \geq C_a(x,\al), \quad a<0.
\ee
Fixing $\al$, we now prove (\ref{eq:ca2}).

Interchanging the order of infimum and supremum in equation (\ref{eq:dpe})
(see \cite{roc}, Corollary 37.3.2), and noting that
the infimum over $m$ is of a continuous function with
compact level sets, we denote by $m^n(y)$ a point where
the infimum is achieved.
By Lemma~\ref{lem:bound:m}, the components
$ \bar\lambda_i^n (y) $ and $ \bar\mu_i^n (y) $ of $ m^n (y) $
are all bounded by $b_0 $.
For $u\in U$ and $y\in G^n$,
\be\label{ineq12}
0 \geq \calL^{n,u,m^n(y)}V^n(y) +c+\rho(u,m^n(y)).
\ee
\begin{lemma}\label{lem:xbar2}
Let $n$ be fixed, and let $\al$ and $x_n$ be as above.
Then there is a filtered probability space
$(\bar\Om,\bar F,(\bar F_t),\bar P)$, and $\bar F_t$-adapted
RCLL processes
$\bar X^n$, $\bar Y^n$ and $u^n$ such that with $\bar P$-probability one
$u^n(t)=\al[\bar m^n](t)$ a.e.\ $t$, $\bar m^n(t)=m^n(\bar X^n(t))$,
$\bar X^n=\Gamma(\bar Y^n)$,
$\bar X^n(0)=\bar Y^n(0)=x_n$,
and for any $f$,
\[
f(\bar X^n(t)) - \int_0^t \calL^{n, u^n(s),m^n}f(\bar X^n(s))ds
\]
\[
f(\bar Y^n(t)) - \int_0^t \calL_0^{n, u^n(s), m^n}f(\bar Y^n(s))ds
\]
are $(\bar F_t)$-martingales.
\end{lemma}

\proof See Section \ref{sec:lemmas}.

Let $\bar m^n(t)= m^n(\bar X^n(t))$
and let $\bar\sigma^n$ be the first
exit time of $\bar X^n$ from $G$.
By (\ref{ineq12}) and Lemma \ref{lem:xbar2}, for any bounded stopping
time $S\le\bar\sig^n$,
\be\label{ineq22}
V^n(x_n) \geq \bar E_{x_n} \left[ V^n(\bar X^n(S))
+\int_0^{S}[c+\rho(\al[\bar m^n](s),\bar m^n(s))]ds \right].
\ee
Denoting $\bar m^n(t)=((\bar\la_i^n(t),\bar\mu_i^n(t))$,
define $\phi^n$ as $\phi^n=\Gamma(\psi^n)$, where
\[
   \psi^n = x+\int_0^\cdot v(\bar u^n,\bar m^n)ds,
\]
and let
\[
\hat\sig^n_a \doteq \inf\{t:\phi^n (t)\not \in G_a\}.
\]
Then the definition of $C_a(x,\alpha)$ implies 
\[
\int_0^{\hat\sig^n_a}[c+\rho(\al[\bar m^n](s),\bar m^n(s))]ds
 \geq C_a(x,\al).
\]
Apply (\ref{ineq22}) with $S=\hat\sig_a^n\w\bar\sig^n\w T$,
with large enough $T$,
using the fact that $V^n\ge0$ to get
\beaa
V^n(x_n) &\geq& \bar E_{x_n} \left[
 \int_0^{\hat\sig^n_a\w\bar\s^n\w T}
 [c+\rho(\al[\bar m^n](s), \bar m^n(s))]ds \right]\\
 &\ge&
\bar E_{x_n}\left[1_{\hat\sig^n_a\le\bar\sig^n}
 \int_0^{\hat\sig^n_a\w T}
 [c+\rho(\al[\bar m^n](s), \bar m^n(s))]ds \right]\\
 &\ge&
\bar P_{x_n}(\hat\s^n_a\le\bar\s^n)C_a(x,\al).
\eeaa
The proof that
$\bar P_{x_n}(\hat\sig^n_a \le \bar\sigma^n)$ tends to one
is analogous to the proof of the that
$\bar P_{x_n}(\hat\sig^n_a \le \bar\sigma^n)\to0$ in the
upper bound. It is therefore omitted. Hence
\[
\liminf_{n \goto \infty} V^n(x_n) \geq C_a(x,\al).
\]
This gives (\ref{eq:ca2}),
and the proof of (\ref{eq:lb}) is established.
\qed

In fact, the value of the game is independent of $b$ for large $b$,
so that the game has a value with the unbounded action space $M$.
As the result depends on Theorems \ref{th:limit}, \ref{th:unique},
\ref{th:solve} we postpone the proof to Section~\ref{sec:lemmas}.
\begin{theorem}\label{th:b}
For all $b\in[b^*,\infty]$, $V^{b,+}=V^+=V^{b,-}=V^-$.
\end{theorem}
\proof See Section \ref{sec:lemmas}.

\section{The PDE}\label{sec:pde}

In this section we show that the upper and lower values of the game
are the unique Lipschitz viscosity solutions of the PDE (\ref{eq:pde}).
Throughout, the parameter $b\in[b^*,\infty)$ is fixed.
Let
\begin{equation}\label{def:H}
 H(q)=\inf_m\sup_u[\lan q,v(u,m)\ran+\rho(u,m)+c].
\end{equation}
It will be useful to
note that the infimum is over the compact set $M^b$,
and the map $(q,u,m)\mapsto [\lan q,v(u,m)\ran+\rho(u,m)+c]$
is continuous.
The PDE of interest is
\be\label{eq:pde}
\lt\{
\begin{array}{ll}
 H(DV(x))=0, & x\in G^o,\\ 
 \lan DV(x),\gamma_i\ran=0, & i\in I(x),\ x\in\pl_+G ,\\ 
 V(x)=0, & x\in\OBdry G.
\end{array}\rt.
\ee
Here, $\gamma_i$ are the directions of constraint that were introduced in Section 2.
\begin{definition}
Let a Lipschitz continuous function $u:X \to\R$ be given
(where $X\subset G$).
We say that $u$ is a subsolution [respectively, supersolution]
to (\ref{eq:pde}) on $X$ if the following conditions hold. Let
$\theta: X \to\R$ be continuously differentiable on $\bar X$.
Let $y\in X$ be a local maximum [minimum] of the map
$x\mapsto u(x)-\theta(x)$. Then
\be\label{eq:subs}
 H(D\theta(y))\vee\max_{i\in I(y)}\lan D\theta(y),\gamma_i\ran\ge0,
\ee
\be\label{eq:supers}
 [\quad
 H(D\theta(y))\w\min_{i\in I(y)}\lan D\theta(y),\gamma_i\ran\le0,
 \quad]
\ee
and
\be\label{eq:subsb}
 V(x)\le0,\quad x\in\bar X\cap\OBdry G,
\ee
\be\label{eq:supersb}
 [\quad
 V(x)\ge0,\quad x\in\bar X\cap\OBdry G.
 \quad]
\ee
We say that $V$ is a viscosity solution to (\ref{eq:pde}),
if it is both a subsolution on $G$ and a supersolution
on $G\setminus\pl_cG$.
\end{definition}

\remark
In case that $\pl_cG\ne\emptyset$, a viscosity solution is often
called a constrained viscosity solution (cf.\ Soner \cite{son},
Capuzzo-Dolcetta and Lions~\cite{capLion}).
The requirement that $V$ is a subsolution up to the boundary
$\pl_cG$---the part of the boundary where exit can be unilaterally
blocked---serves as a boundary condition on this
part of the boundary. Note that in the current paper,
the term `constrained' refers to the part $\pl_+G$ of the
boundary, where it is the mechanism associated with the Skorokhod
Problem that constrains the dynamics to $G$.
\qed

First, we address uniqueness of solutions to (\ref{eq:pde}).
\begin{theorem}\label{th:unique}
Let $u$ be a subsolution and $v$ a supersolution to (\ref{eq:pde}).
Then $u\le v$ on $G$.
\end{theorem}

The proof combines ideas from two sources, namely \cite{atadup2} (which is based on \cite{dupish2}, and discusses how to deal  with the constrained
dynamics on $\pl_+G$), and \cite{son} (to accommodate the fact that under Condition \ref{cond:G}.1
part of the boundary ($\pl_cG$) can be thought
of as imposing a state-space constraint on the maximizing player).

The following lemma will be used in proving Theorem \ref{th:unique}.
In the interest of consistency with previous
publications, we use $B$ in Lemma~\ref{lem:mu} below to denote a certain
subset of $\R^J$ (although everywhere except in this section, $B$ denotes
a set of strategies).
Part 1 states that the ``Set B'' Condition holds, namely a condition
under which it was proved in \cite{dupish1} that the SM enjoys
the regularity property (\ref{sp:lip}).
The proof that this condition holds in
the current setting can be found in \cite{dupram23}.
The existence of a smooth version of the set $B$ is proved in
\cite{atadup2}
(before Lemma 2.1). For Parts 2 and 3, see Lemmas 2.1 and 2.2 of
\cite{atadup2} (note that the condition that $\gamma_i$ are independent
holds).
\begin{lemma}\label{lem:mu}
\ben
\item
There exists a compact, convex, and symmetric set $B\subset \RJ $ with
$0\in B^o$, such that if $z\in\pl B$ and if $n$ is an outward
normal to $B$ at $z$, then for all $i\in\{1,\ldots,J\}$
$$
 \lan z,e_i\ran\ge-1 {\rm\ implies\ } \lan\gamma_i,n\ran\ge0
 {\rm\ and\ } \lan z,e_i\ran\le1 {\rm\ implies\ } \lan\gamma_i,n\ran\le0.
$$
In addition, the unit outward normal $n(x)$ to $B$ at $x\in\pl B$ is
unique and continuous (as a function on $\pl B$).
\item
Let $\bar n$ be the extension of $n$ to $ \RJ $
satisfying $\bar n(x)=n(y)$ whenever $ax=y\in\pl B$, some $a\in(0,\infty)$
(and define $\bar n(0)$ arbitrarily).
Let $\Xi: \RJ \to\R_+$ be defined via
$$
 \Xi(x)=a\ \Leftrightarrow\ x\in\pl(aB)
$$
for all $a\in[0,\infty)$, and let $\xi(x)=(\Xi(x))^2$.
Then there exist constants $m,M\in(0,\infty)$ and a function
$\vr: \RJ \to[m,M]$ such that the $C^1( \RJ )$
function $\xi$ satisfies $m\|x\|^2\le\xi(x)\le M\|x\|^2$, and
$D\xi(x)=\vr(x)\Xi(x)\bar n(x)$.
\item
There exists a constant $m_1\in (0,\infty)$ and a continuously differentiable
function $\mu: \RJP \to[0,m_1]$ such that $\|D\mu\|\le m_1$ on
$ \RJP $, and
$$
 \lan D\mu(x),\gamma_i\ran<0,\quad x\in \RJP ,\quad i\in I(x).
$$
\een
\end{lemma}

In what follows, we keep the notation of Lemma \ref{lem:mu}
for $B, \bar n, \Xi, \xi, \vr$ and $\mu$.

\noindent{\bf Proof of Theorem \ref{th:unique}:}
For $a>0$, let
$$
 U(x)=u(x)-a\mu(x),
$$
$$
 V(x)=v(x)+a\mu(x).
$$
Let $\del>0$. Then it suffices to show that for all small $a>0$,
$\del>0$, one has $U\le (1+\del)V$ on $G$. Arguing by contradiction,
we assume that this is not true. Then there are $a$ and $\del$
arbitrarily small such that 
$$
 \rho=\sup_{x\in G}[U(x)-(1+\del)V(x)]>0.
$$
Below we let $c_i, i=1,2,...$ denote positive constants.
Consider  Condition \ref{cond:G}.1 first.
Let
\begin{equation}\label{eq:phixy}
 \Phi(x,y)= U(x)-(1+\del)V(y)-\frac 1 \eps\xi(x -y-\eps^{1/2}y).
\end{equation}
Let $(\bar x,\bar y)\in \bar G^2$
achieve the maximum of $\Phi$ in $ \bar G \times
\bar G $.
By continuity of $U$ and $V$, there exists $ \bar z\in\bar G $ so that 
$ \rho = U( \bar z)-(1+\del)V(\bar z) $. Note that
\begin{equation}\label{eq:10}
(1+\eps^{1/2})^{-1}\bar z\in\bar G.
\end{equation}
Hence by the Lipschitz continuity
of $V$,
\bea\label{e:phi1}\nonumber
\Phi ( \bar x , \bar y )
& \geq & \Phi 
         \left ( \bar z, \frac {\bar z} {1 + \eps^{\Sfrac 1 2 }} \right )
\\
\nonumber
& =  & U(\bar z)-(1+\del)
V \left (\frac {\bar z} {1 + \eps^{\Sfrac 1 2 }} \right )
\\
& \geq & \rho - c_1 \eps^{\Sfrac 1 2 } .
\eea
By Lipschitz continuity of $U$ and the lower bound on $\xi$ given in Lemma \ref{lem:mu}, 
\bea\label{e:phi2}\nonumber
\Phi ( \bar x , \bar y ) 
& =    & U(\bar x)-(1+\del)V(\bar y)
         -\frac 1 \eps \xi (\bar x  - \bar y-\eps^{1/2}\bar y)   \\
\nonumber
& \leq & U (\bar y ) + c_2 | \bar x - \bar y | - (1+\del)V(\bar y)
         - \frac m \eps | \bar x - \bar y-\eps^{1/2}\bar y |^2\\
&\le&
\rho + c_2 | \bar x - \bar y | 
     - \frac{c_3}{\eps} | \bar x- \bar y-\eps^{1/2}\bar y |^2 .
\eea
By (\ref{e:phi1}) and (\ref{e:phi2}),
\bea\label{e:phi3}
c_2 | \bar x - \bar y | + c_1 \eps^{\Sfrac 1 2}
&\geq&
\frac{c_3}{\eps} | \bar x- \bar y-\eps^{1/2}\bar y |^2
\\ \nonumber
&\geq&
\frac{c_4}{\eps} | \bar x - \bar y |^2 - c_4 | \bar y |^2 .
\eea
Since $ \bar x $ and $ \bar y $ are bounded,
(\ref{e:phi3}) implies
$| \bar x - \bar y |^2 \leq c_5 \eps$
and so
\begin{equation}\label{e:phi5}
| \bar x - \bar y | \leq c_6 \eps^{\Sfrac 1 2 }.
\end{equation}
Using this in~(\ref{e:phi3}) we have
\begin{equation}\label{e:phi4}
| \bar x- \bar y-\eps^{1/2}\bar y | \leq c_7 \eps^{\Sfrac 3 4} .
\end{equation}
By (\ref{e:phi5}),
$\bar x-\bar y\to0$ as $\eps\downarrow0$. Also,
we claim that
for all $\eps>0$ small, $\bar x$ and $\bar y$ are 
bounded away from $\pl_o G$.
To see this, assume the contrary. Then along a subsequence, both
$\bar x$ and $\bar y$ must converge to the same point on $\pl_o G$.
Using the continuity of $u$ and $v$,
(\ref{eq:subsb})--(\ref{eq:supersb}),
and the non-negativity of $\xi$,
$\limsup\Phi(\bar x,\bar y)
\le \limsup \left [ u(\bar x)-(1+\del)v(\bar y) \right ] \le0$,
where the limit superior is taken along this subsequence.
However, by (\ref{e:phi1}),
for all small $\eps$, $\Phi(\bar x,\bar y)\ge\rho/2>0$, which
gives a contradiction.

Let
$$
 \theta(x)=\frac1\eps\xi(x-\bar y-\eps^{1/2}\bar y)+a\mu(x),
$$
and note that the map $x\mapsto u(x)-\theta(x)$ has a maximum
at $\bar x\in G$.
Since $u$ is a subsolution, (\ref{eq:subs}) must
be satisfied at $\bar x$. Denoting
\begin{equation}\label{eq:que}
q^\eps=
\vr(\bar x-\bar y-\eps^{1/2}\bar y)\Xi(\bar x-\bar y-\eps^{1/2}\bar y)
\bar n(\bar x-\bar y-\eps^{1/2}\bar y),
\end{equation}
we have from Lemma \ref{lem:mu}.2 that
$$
 D\theta(\bar x)=\frac1\eps q^\eps+aD\mu(\bar x).
$$
Suppose $i$ is such that $\bar x_i=0$. Then 
$ \lan\bar x-\bar y-\eps^{1/2}\bar y,e_i\ran\le0 $
and so by Lemma \ref{lem:mu}.1,
$$
\lan\gamma_i, \bar n(\bar x-\bar y-\eps^{1/2}\bar y)\ran\le0.
$$
Since by Lemma \ref{lem:mu}.3, $\lan\gamma_i, D\mu(\bar x)\ran<0$,
it follows that $\lan\gamma_i,D\theta(\bar x)\ran<0$. It follows from
(\ref{eq:subs}) that
$$
 H(D\theta(\bar x))\ge0,
$$
namely,
\be\label{eq:H1}
 H\lt(\frac1\eps q^\eps+aD\mu(\bar x)\rt)\ge0.
\ee
On the other hand, let
$$
 \al(y)=-a\mu(y)-\frac{1}{\eps(1+\del)}\xi(\bar x-y-\eps^{1/2}y).
$$
Note that
$$
 D\al(\bar y)=-aD\mu(\bar y)+\frac{1+\eps^{1/2}}{\eps(1+\del)} q^\eps
$$
and that
the map $y\mapsto v(y)-\al(y)$ has a minimum at $\bar y$.
Since $\bar x\in\bar G$, (\ref{e:phi4}) implies that
$\bar y\in G\setminus\pl_cG$ for all small $\eps$.
Since $v$ is a supersolution,
(\ref{eq:supers}) is satisfied at $y$. An argument as above shows
that
$$
 H(D\al(\bar y))\le0,
$$
and therefore
$$
 H\lt(\frac{1}{1+\del}\lt(\frac{1+\eps^{1/2}}\eps q^\eps
 -a(1+\del)D\mu(\bar y)\rt)\rt)\le0.
$$
It follows from the definition of $H$, using $\rho(u,m)\ge0$, that
$$
 H\lt(\frac 1{1+\del}\, p\rt)\ge\frac 1{1+\del}H(p)+\frac\del{1+\del}c,
$$
and therefore
\be\label{eq:H2}
 H\lt(\frac1\eps q^\eps+\frac1{\eps^{1/2}}q^\eps
-a(1+\del)D\mu(\bar y)\rt)+\del c\le0.
\ee
Now $D\mu$ is bounded, and by (\ref{e:phi4}), boundedness of
$n$ and $\vr$, and
the Lipschitz continuity of $\Xi$, it follows that
$\eps^{-1/2}q^\eps$ converges to zero as $\eps\to0$.
Note that by (\ref{def:H}) and the following comment,
$H$ is uniformly continuous on $\R^J$.
Therefore,
(\ref{eq:H1}) and (\ref{eq:H2}) give a contradiction when $a>0$
and $\eps>0$ are small and $\del>0$ fixed.

Under Condition \ref{cond:G}.2 the above argument
is not valid, since (\ref{eq:10}) may not hold.
However,
in this case the minimizing player can force exit from any point on $\partial G\backslash \partial_+ G$,
and the additional complications due to the ``state-space constraint'' used under part 1 are no longer needed.
In other words,
instead of (\ref{eq:phixy})  we can consider 
$$
\Phi(x,y)=U(x)-(1+\del)V(y)-\frac1\eps\xi(x-y),
$$
and a review of the above proof shows that (\ref{eq:H1})
and (\ref{eq:H2}) still hold if the expression
$\eps^{1/2}$ is replaced by zero everywhere in (\ref{eq:que}) and
(\ref{eq:H2}). A contradiction is then obtained analogously.
\qed

We next consider the upper and lower values of the game,
and remind the reader that in this section the rates $m$ are assumed bounded.
\begin{theorem}\label{th:solve}
$V^-$ and $V^+$ are solutions to (\ref{eq:pde}).
\end{theorem}
Recall that from Lemma \ref{lem:vallip}, $V^\pm$ are Lipschitz.

\noindent
{\bf Proof of Theorem \ref{th:solve}:}
We use the specific form of $v(u,m)$ and $\rho(u,m)$. These can be
written as
$$
 v(u,m)=b_0(m)+\sum_{i=1}^Ju_ib_i(m),
$$
$$
 \rho(u,m)=c_0(m)+\sum_{i=1}^Ju_ic_i(m).
$$
We have $\sum_{i\in C(k)}u_i\le1$.
The $b_i$ are linear, and the $c_i$ are convex.
Hence, as a direct consequence of \cite[Corollary 37.3.2]{roc}, 
the {\em Isaacs condition} holds, namely
\be\label{eq:isaacs}
 H(q)=\inf_m\sup_u[\lan q,v(u,m)\ran+c+\rho(u,m)]
 =\sup_u\inf_m[\lan q,v(u,m)\ran+c+\rho(u,m)].
\ee
Another fact that we will use is that for any
$y\in G\setminus\pl_cG$ there
is $\del_0=\del_0(y)>0$ which serves as a lower bound on the exit time.
Namely, if $\phi$ solves $\dot\phi=\pi(\phi,v(u,m))$, $\phi(0)=y$,
then
\be\label{eq:boundsig}
 \sig \doteq \inf\{ t\geq 0: \phi(t) \not \in G\} \ge\del_0,\quad u\in \bar U, m\in \bar M.
\ee
The bound is an immediate consequence of the $u$ and $m$ being
uniformly bounded.

By definition, $ V ^\pm (x) = 0 $ for $ x \in \pl_o G $.
Thus we only need to
establish (\ref{eq:subs})--(\ref{eq:supers}).
The proof consists of four parts.

\noindent{\em \underline{Proof that $V^-$ is a supersolution on $G \backslash \partial_cG$.}}

Standard dynamic programming arguments show that for $\del>0$,
\be\label{eq:dp1}
 V^-(x)=\inf_\beta\sup_u\lt[\int_0^{\s\w\del}(c+\rho(u,\beta[u]))dt
 +V^-(\phi(\s\w\del))\rt],
\ee
where $\phi$ is the solution to $\dot\phi=\pi(\phi,v(u,\beta[u]))$,
with $\phi(0)=x$.
Let $\theta$ be smooth, and let $y\in G\setminus\pl_cG$
be a local minimum of
$V^--\theta$. We can assume without loss that $V^-(y)=\theta(y)$.
We need to show
\be\label{alpha}
 H(D\theta(y))\w\min_{i\in I(y)}\lan D\theta(y),\gamma_i\ran\le0.
\ee
We shall assume the contrary and reach a contradiction.
Thus, there exists $a>0$ such that
$H(D\theta(y))\ge a$, and
\be\label{eq:D1}
 \lan D\theta(y),\gamma_i\ran\ge a,\quad i\in I(y).
\ee
 From the definition of $H$ and (\ref{eq:isaacs}),
$$
 \sup_u\inf_m[\lan D\theta(y),v(u,m)\ran+c+\rho(u,m)]\ge a,
$$
and therefore there exists a $u_0$ such that for all $m$,
$$
 \lan D\theta(y),v(u_0,m)\ran+c+\rho(u_0,m)\ge a/2.
$$
For any strategy $\beta$, if $\bar u(t)\equiv u_0$,
\be\label{eq:D2}
 \lan D\theta(y),v(\bar u(t),\beta[\bar u](t))\ran+c
 +\rho(\bar u(t),\beta[\bar u](t))\ge a/2
\ee
for all $t$. Let $\phi$ denote the dynamics corresponding to $\bar u$
and a generic $\beta$, starting from $y$.
Note that the mapping $z\mapsto I(z)$ is upper semi-continuous, in
the sense that for any $z$ there is a neighborhood of $z$ on which
$I(\cdot)\subset I(z)$. Using the boundedness of $m$
this implies that for any $\beta\in B$, one has
$I(\phi(r))\subset I(y)$ for $r\in[0,\del]$, if $\del>0$ is
chosen small enough.
We now use that $\phi$ is a solution to the SP.
Choosing such a $\del >0$, for any $r\in[0,\del]$ there
exist $a_i\ge0$ (that may depend on $r$) such that
$$
 \dot\phi(r)=v(\bar u(r),\beta[\bar u](r))+\sum_{i\in I(y)}a_i\gamma_i.
$$
Using the continuity of $D\theta$ and taking $\del>0$ smaller
if necessary, (\ref{eq:D1}) and (\ref{eq:D2}) imply,
for $t\in[0,\del]$,
\beaa
 \frac{d}{dt}\theta(\phi(t)) &=&
\lan D\theta(\phi(t)),\dot\phi(t)\ran\\
 &\ge&
-c-\rho(\bar u(t),\beta[\bar u](t))+a/4.
\eeaa
Taking $\del$ even smaller if necessary (so that it is at most
$\del_0$), we have from (\ref{eq:boundsig}) that
$$
 \theta(\phi(\del))-\theta(y)\ge
 -\int_0^{\del}(c+\rho(\bar u(t),\beta[\bar u](t)))dt
 +a\del/4.
$$
 From (\ref{eq:dp1}), one can find a $\beta$ such that
$$
 V^-(y)\ge\sup_u\lt[\int_0^{\del}(c+\rho(u,\beta[u]))dt
 +V^-(\phi(\del))-a\del/8\rt].
$$
Letting $u=\bar u$, the last two displays give (using $\theta(y)
=V^-(y)$)
$$
 \theta(\phi(\del))\ge V^-(\phi(\del))+a\del/8,
$$
so that $V^-(\phi(\del))-\theta(\phi(\del))<0$ for all
$\del>0$ small, contradicting the assumption that $y$ is
a local minimum of $V^--\theta$. This proves that $V^-$ is a
supersolution on $G \backslash \partial_cG$.

\noindent
{\em \underline{Proof that $V^-$ is a subsolution on $G$.}}

Let $\theta$ be smooth and $y\in G$ a local maximum of
$V^--\theta$.
In case that $y\in\pl_cG$, let $\bar U_{y,\beta,\del}$
be the set
of controls $u\in\bar U$ for which the trajectory $\phi$ determined
by $u$ and $\beta[u]$ and starting from $y$ does not exit $G$
on $[0,\del]$. Given $y\in\pl_cG$, it is clear that
$\bar U_{y,\beta,\del}$ is not empty for all $\del$
small and all $\beta$, by considering the control $u=0$.
Moreover, for all $\del$ small enough,
(\ref{eq:dp1}) is valid where the supremum
extends only over $u\in\bar U_{y,\beta,\del}$. Indeed, given
$u\not\in\bar U_{y,\beta,\del}$, consider
$u'$ that agrees with $u$ on $[0,\sig]$ and $u'=0$ on
$(\sig,\del]$. Then the expression in brackets in (\ref{eq:dp1})
is identical under $u$ and under $u'$, but
$u'\in\bar U_{y,\beta,\del}$.

Assume without loss that $V^-(y)=\theta(y)$.
We would like to show that
\be\label{gamma}
 H(D\theta(y))\vee\max_{i\in I(y)}\lan D\theta(y),\gamma_i\ran\ge0.
\ee
Assuming the contrary, there exists $a>0$ such that
$H(D\theta(y))\le-a$, and
\be\label{eq:D3}
 \lan D\theta(y),\gamma_i\ran\le-a,\quad i\in I(y).
\ee
Using the definition of $H$ and (\ref{eq:isaacs}), for all $u$ there exists $m_u$ such that
\be\label{eq:D4}
\lan D\theta(y),v(u,m_u)\ran
 +c+\rho(u,m_u)\le-a/2.
\ee
Note that it is possible to choose $m_u$ so that it depends continuously
on $u$.
Define $\bar\beta$ as $\bar\beta[u](t)=m_{u(t)}$ for all $t$. Since
$\bar\beta[u]$ is measurable if $u$ is, $\bar\beta$ maps $\bar U$
into $\bar M$.
Let $\phi$ be the trajectory corresponding to $\bar\beta$ and a generic
$u\in \bar U$, (or a generic $u\in\bar U_{y,\beta,\del}$ if
$y\in\pl_cG$) starting from $y$.
Arguing as before by upper semi-continuity
of $I(\cdot)$, if $\del$ is small enough, then
$$
 \dot\phi(r)=v(u(r),\bar\beta[u](r))+\sum_{i\in I(y)}a_i\gamma_i,
 \quad r\in[0,\del],
$$
where $a_i\ge0$ may depend on $r$.
By possibly taking $\del$ smaller, and smaller
than $\del_0$, we have, using the continuity of $D\theta$ and
(\ref{eq:D3}), (\ref{eq:D4}) that
\beaa
 \frac{d}{dt}\theta(\phi(t)) &=&
\lan D\theta(\phi(t)),\dot\phi(t)\ran\\
 &\le&
-c-\rho(u(t),\bar\beta[u](t))-a/4,
\eeaa
and
$$
 \theta(\phi(\del))-\theta(y)\le-\int_0^\del(c+
 \rho(u(t),\bar\beta[u](t)))dt-a\del/4.
$$
Now, (\ref{eq:dp1}) implies that for any $\beta$ there
is $u$ such that
$$
 V^-(y)\le\int_0^\del(c+\rho(u,\beta[u]))dt+V^-(\phi(\del))+a\del/8.
$$
Specializing to $\beta=\bar\beta$, the last two displays show
that $V^-(\phi(\del))-\theta(\phi(\del))>0$ for all $\del>0$
small. This contradicts the assumption that $y$ is a local
maximum of $V^--\theta$, and as a result, $V^-$ is a subsolution.

\noindent
{\em \underline{Proof that $V^+$ is a supersolution on $G \backslash \partial_cG$.}}

The proof is analogous to the proof that $V^-$ is a subsolution.
Most details are therefore skipped.
The dynamic programming principle states that
for $\del>0$,
\be\label{eq:dp2}
 V^+(x)=\sup_\al\inf_m\lt[\int_0^{\s\w\del}(c+\rho(\al[m],m))dt
 +V^+(\phi(\s\w\del))\rt],
\ee
where $\phi$ is the dynamics corresponding to $\al$ and $m$,
starting from $x$. Taking a smooth $\theta$, and leting
$y\in G\setminus\pl_cG$
be a local minimum of $V^+-\theta$, showing
$$
 H(D\theta(y))\w\min_{i\in I(y)}\lan D\theta(y),\gamma_i\ran\le0
$$
can be obtained by an argument analogous to that used to prove (\ref{gamma}),
using (\ref{eq:dp2}) in place of (\ref{eq:dp1}).

\noindent
{\em \underline{Proof that $V^+$ is a subsolution on $G$.}}

We need to show that
\be\label{beta}
 H(D\theta(y))\vee\max_{i\in I(y)}\lan D\theta(y),\gamma_i\ran\ge0,
\ee
where $\theta$ is smooth, and $y\in G$ is a local maximum of
$V^+-\theta$.
In the special case where $y\in\pl_cG$, we can assume
without loss that the supremum in (\ref{eq:dp2}) extends
only over $\al\in A_{y,\del}$, the set of strategies
under which, for any $m\in\bar M^b$, the dynamics associated
with $\al$ and $m$, and starting from $y$,
does not leave $G$ before $\del$.
The proof of (\ref{beta})
is analogous to the proof of (\ref{alpha}), and is skipped.

This completes the proof that $V^-$ and $V^+$ are solutions
to (\ref{eq:pde}).
\qed

\section{A competing queues example}\label{sec:example}

Consider a queueing network with only one server, providing
service to $J$ classes. Each customer requires service once.
In this example all arrival rates are positive:
$\la_i>0$ for all $i$, hence $\calJ_+=\{1,\ldots,J\}$.
This network, ``the $k$ competing queues,'' has been studied extensively,
in discrete and continuous
time (see~\cite{BDM,W} and references therein). When the criterion 
(to be minimized) is
either the average cost or the discounted cost, and the one-step cost
is a positive linear combination
$\sum_i c_ix_i$ of the queue sizes $x_i$, the optimal policy is
the $ \mu $-$c$ rule, which is a priority discipline,
giving absolute priority to the non-empty queue for which
$\mu_ic_i$ is maximal.
Under the cost studied here, the optimal policy is quite different.

\begin{figure} 
\centerline{
\begin{tabular}{cc}
\psfig{file=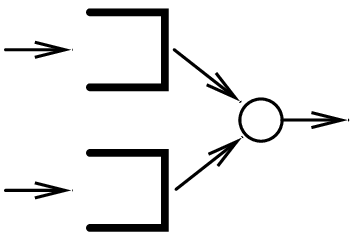}\qquad\qquad\qquad
\psfig{file=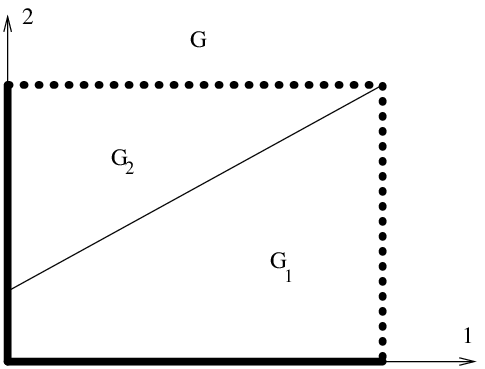}
\end{tabular}
}
\caption{Priority to class $i$ when the state is in $G_i$, $i=1,2$.}
\label{fig:example}
\end{figure}

\begin{proposition}
Consider the case where $G$ is a hyper-rectangle, given
as $G=\{x:0\le x<z_i\}$, where $z_i>0$ are constants.
Assume that $\la_i>0$ for all $i=1,\ldots,J$.
If $c$ is large enough,
then the viscosity solution to the PDE (\ref{eq:pde}) is given as
\begin{equation}\label{eq:form}
V(x)=\min_i\al_i(z_i-x_i),
\end{equation}
where $\al_i>0$ are constants depending on $c$.
\end{proposition}
We remark that the constants $\al_i$ are uniquely defined by
(\ref{eq:form1}) below.
In the totally symmetric case, where $\mu_i=\mu$, $\la_i=\la$,
$z_i=z$ for all $i$, the solution takes the form
$V(x)=\al\min_i(z-x_i)$. In this case, the optimal
service discipline can be interpreted as ``serve the longest
queue.'' An asymmetric two dimensional example
is given in Figure \ref{fig:example}, where the domain $G$
is divided into two subdomains $G_1$ and $G_2$ in accordance
with the structure (\ref{eq:form}), and the optimal service
discipline corresponds to giving priority to
class $i$ when the state is within $G_i$, $i=1,2$.
Thus the optimal control under our escape-time criterion is very different
from the optimal controls for the average or discounted cost criteria.

\proof
The constraint directions are given by $\gamma_i=e_i$.
The Hamiltonian is given by
$$
H(p)=\sup_u\inf_m H(p,u,m),
$$
where
$$
H(p,u,m)=c+\sum_i\lt[p_i(\bar\la_i-u_i\bar\mu_i)
+\la_i l\left (\frac{\bar\la_i}{\la_i}\right)
+u_i\mu_il\left(\frac{\bar\mu_i}{\mu_i}\right)\rt].
$$
Using strict convexity and smoothness of the map
$m\mapsto H(p,u,m)$, the minimum over $m$ is attained
at $\bar\la_i=\la_i e^{-p_i}$, $\bar\mu_i=\mu_ie^{p_i}$.
Thus
$$
H(p,u)\doteq \inf_m H(p,u,m)=
c+\sum_i[\la_i(1-e^{-p_i})+u_i\mu_i(1-e^{p_i})].
$$
For the proposed solution, $DV(x)\in-\R_+^J$
wherever the gradient is defined. For $p\in-\R_+^J$, maximizing
$H(p,u)$ over $u$ clearly gives
\begin{equation}\label{eq:ham}
H(p)=\sup_u H(p,u)=c+\sum_i\la_i(1-e^{-p_i})+\max_i\mu_i(1-e^{p_i}).
\end{equation}
We use the well known fact that the definition
of viscosity solutions can be equivalently stated in terms
of sub- and superdifferentials (see \cite{barcap}, Lemma II.1.7).
Note that (\ref{eq:subsb}) and (\ref{eq:supersb})
hold, since $V=0$ on $\pl_oG$.
Hence it suffices to verify
that (\ref{eq:subs}) [resp., (\ref{eq:supers})] holds where
$D\theta(y)$ is replaced by any superdifferential [subdifferential]
of $V$ at $y$. 

We show first that the equation $H(DV(x))=0$ holds wherever
$DV$ is defined. 
The proposed form (\ref{eq:form}) satisfies $DV(x)=-\al_ie_i$,
wherever the gradient is defined, with $i=i_x$ depending on $x$.
By the special form of the gradient,
the equation $H(DV(x))=0$ takes the form
\begin{equation}\label{eq:form1}
H(DV(x))=c+\la_i(1-e^{\al_i})+\mu_i(1-e^{-\al_i})=0,
\end{equation}
where $i=i_x$. Denote $c_i=c/(\la_i+\mu_i)$.
Then equivalently,
$
1+c_i-F_i(\al_i)=0,
$
where
$$
F_i(\al_i)=\frac{\la_ie^{\al_i}+\mu_ie^{-\al_i}}{\la_i+\mu_i}.
$$
The function $F_i$ is strictly convex, $F_i(0)=1$,
and $F_i(\al_i)\to\infty$ as $\al_i\to\infty$. Since $c_i>0$,
it follows that there are unique positive constants
$\al_i$ where $F_i(\al_i)=1+c_i$, $i=1,\ldots,J$.
These are the constants 
in (\ref{eq:form}). In particular, (\ref{eq:form1}) holds
for $i=i_x$, and $H(DV(x))=0$.

Next consider any interior point $x$ at which the
gradient is not defined. Clearly
there are no subdifferentials at that point, and any
superdifferential is given as a convex combination of
$-\al_ie_i$, $i=1,\ldots,J$.
Let $B(e_i,\eps)$ be the open ball of radius $\eps$ about $e_i$.
Denote
$\til S=\{\nu\in\R^J:\nu_i\ge0,\sum\nu_i=1\}$,
$S=\{\nu\in\R^J:\nu_i\ge0,\sum\nu_i\le1\}$,
$S_\eps=S\cap\cup_iB(e_i,\eps)$, and $S_\eps^c=S-S_\eps$.
Let $q=-\sum_i\nu_i\al_ie_i$. It suffices to show
that $H(q)\ge0$ for $\nu\in\til S$, but since we later need a stronger
statement than that, we show that in fact $H(q)\ge0$
holds for $\nu\in S$. By (\ref{eq:ham}),
$$
H(q)=c+\sum_i\la_i(1-e^{\nu_i\al_i})+\max_i\mu_i(1-e^{-\nu_i\al_i}).
$$
Define
$$
H^1(q)=c+\sum_i\la_i(1-e^{\nu_i\al_i})+\mu_1(1-e^{-\nu_1\al_1}),
$$
$$
\bar H(q)=c+\sum_i[\la_i(1-e^{\nu_i\al_i})+\mu_i(1-e^{-\nu_i\al_i})].
$$
By (\ref{eq:form1}), $c+\la_i(1-e^{\al_i})+\mu_i\ge0$ and
$c+\la_i(1-e^{\al_i})\le0$, and it follows that there are
constants $A_1,A_2,A_3$ and $A_4$ such that for all $c$ and $i=1,\ldots,J$,
\begin{equation}\label{eq:A}
A_1+\log(c+A_2)\le\al_i\le A_3+\log(c+A_4).
\end{equation}

We first consider small perturbations $\nu$ of $e_1$.
To show that $H(q)\ge0$, it suffices to show that $H^1(q)\ge0$.
Note that (\ref{eq:form1}) implies $H^1(q)|_{\nu=e_1}=0$. Also,
$$
\nabla_\nu H^1(q)|_{\nu=e_1}
=(-\la_1\al_1e^{\al_1}+\mu_1\al_1e^{-\al_1})e_1
-\sum_{i\ne1}\la_i\al_ie_i.
$$
Hence, for $\gamma=e_i-e_1$
(where $i\ne1$), using $c+\la_1(1-e^{\al_1})\le0$,
(\ref{eq:form1}) and (\ref{eq:A}),
\beaa
\nabla_\nu H^1(q)|_{\nu=e_1}\cdot\gamma
&=&
\la_1\al_1e^{\al_1}-\mu_1\al_1e^{-\al_1}-\al_i\la_i\\
&=&
\al_1(2\la_1e^{\al_1}-\mu_1-c-\la_1)-\al_i\la_i\\
&\ge&
\al_1(c-\mu_1)-\al_i\la_i\\
&\ge&
[A_1+\log(c+A_2)](c-\mu_1)-[A_3+\log(c+A_4)]\la_i\\
&\ge&
1,
\eeaa
for all $c$ large. 
Analogous calculations give $\nabla_\nu H^1(q)|_{\nu=e_1}
\cdot\gamma\ge1$ for $\gamma=-e_1$ as well.
As a result, the directional derivatives
$(\pl/\pl\til\gamma)H^1(q)|_{\nu=e_1}$ in the direction $\til\gamma$,
where $\til\gamma$ are of the form
$\til\gamma=(y-e_1)/\|y-e_1\|$, $y\in S$, are bounded below
by $1/2$. Hence $H^1(q)\ge0$ for $\nu\in S$ within
a neighborhood of $e_1$
and $c$ large. Consequently, a similar statement holds for
$H(q)$. Since the same argument holds for neighborhoods of
$e_i$, $i=2,\ldots,J$, we conclude that there is $\eps>0$
and $c_0$ such that $H(q)\ge0$ for $\nu\in S_\eps$ and $c\ge c_0$.

Next consider $\nu\in S_\eps^c$.
We first provide a lower bound on $(\pl/\pl c)\bar H(q)$.
Differentiating (\ref{eq:form1}) with respect to $c$,
$\dot\al_i\doteq
\pl\al_i/\pl c=(\la_i e^{\al_i}-\mu_ie^{-\al_i})^{-1}$.
Using (\ref{eq:A}), for all $c$ large,
$0\le\dot\al_i\le(\la_ie^{\al_i}-1)^{-1}$. Using this,
the fact that $\nu$ is bounded away from $\cup_i\{e_i\}$,
and by taking $c$ large, one has
\beaa
\frac{\pl}{\pl c}\bar H(q)
&\ge&
1-\sum_i\nu_i\la_i\dot\al_ie^{\nu_i\al_i}\\
&\ge&
1-\sum_i\nu_i[e^{(1-\nu_i)\al_i}-1]^{-1}\\
&\ge&
\frac12.
\eeaa
Note that the above bound holds for all $c\ge c_1$ and
all $\nu\in S_\eps^c$, where $c_1$ is a constant.
It follows that there is $c_2$ such that for all $c\ge c_2$
and all $\nu\in S_\eps^c$, one has $\bar H(q)\ge\sum_i\mu_i$.
Since $H\ge \bar H-\sum_i\mu_i$, $H(q)\ge0$.
We conclude that $H(q)\ge0$ for all $\nu\in S$.
In particular, $H(q)\ge0$ where $q$ is any
superdifferential of $V$ at any interior point.

Finally, consider a point $x\in G\cap\pl\R_+^J$. Any superdifferential
of $V$ at $x$ is given as
$q=\sum_{i\in I(x)}\eta_ie_i-\sum_{j=1}^J\nu_j\al_je_j$,
where $\eta_i\ge0$.
If $\max_{i\in I(x)}\lan q,\gamma_i\ran\ge0$, then
(\ref{eq:subs}) holds. Otherwise, $\lan q,\gamma_i\ran<0$
for all $i\in I(x)$. Consequently, any $q$ of the form above
is given as $-\sum_{j=1}^J\nu_j'\al_je_j$, with $\nu'\in S$.
As we have shown, in this case, $H(q)\ge0$. Therefore
(\ref{eq:subs}) holds.

Similarly, any subdifferential of $V$ at $x\in G\cap\pl\R_+^J$ is
of the form $-\sum_{i\in I(x)}\eta_ie_i-\sum_{j=1}^J\nu_j\al_je_j$.
In particular, $\lan q,e_i\ran\le0$ for all $i$, and (\ref{eq:supers})
holds.
\qed

\section{Proofs of lemmas}\label{sec:lemmas}

\noi{\bf Proof of Lemma \ref{lem:dpe}:}
Let 
\[
W^{n}(x)\doteq \inf E_{x}^{u,n}e^{-nc\sigma _{n}}. 
\]
Since $c>0$, $W^{n}$ is well defined.  
Standard iterative methods can be used to construct a solution to the DPE 
\begin{equation}\label{eq:dpe3}
0=\inf_{u\in U}\left[ \til {\cal L}^{n,u}\bar W^{n}(x)-nc\bar W^{n}(x)\right] ,x\in G^{n}
\end{equation}
and the boundary condition $\bar W^{n}(x)=1$ if $x\notin G^{n}$.
We claim that this solution coincides with the risk-sensitive cost.
 To see this,
consider a controlled Markov process $(X^n,u)$ that starts at $x$.
Then
\[
Y(t) \doteq \bar W^n(X^n(t))-\bar W^n(x)
 -\int_0^t\tilde{\cal L}^{n,u(s)}\bar W^{n}(X^n(s))ds
\]
is a martingale.
Equation (\ref{eq:dpe3}) implies 
${\cal L}^{n,u(s)}\bar W^{n}(X^n(s))\geq nc\bar W^{n}(X^n(s))$,
and so 
\[
\bar W^n(X^n(t))-\bar W^n(x)-\int_0^t nc \bar W^{n}(X^n(s))ds
=\int_0^tZ(s)ds+Y(t)
\]
for some nonnegative process $Z$.
Using Gronwall's lemma we obtain that for each $t<\infty$
\[
E^{n,u}_x \bar W^n(X^n(t\wedge \sigma^n))e^{-nc(t\wedge \sigma^n)}\geq \bar W^n(x),
\]
and by the Lebesgue Dominated Convergence Theorem 
\[
E^{n,u}_x e^{-nc\sigma^n}\geq \bar W^n(x).
\]
If we define $u$ in terms of the feedback control that minimizes in (\ref{eq:dpe3}) then all the inequalities above become equalities, thus showing that $\bar W^n=W^n$.

The definition of $W^{n\,}$ implies $W^{n}(x)=\exp \left[ -nV^{n}(x)\right] $
. \ If we insert this into the DPE of $\bar W^{n}$ and multiply by $\exp \left[
nV^{n}(x)\right] $ then the equation 
\begin{eqnarray*}
0& =& \inf_{u\in U}\left[ \sum_{j=1}^{J}n\lambda _{j}\left( \exp \left[
-nV^{n}\left( x+\frac{1}{n}v_{j}\right) +nV^{n}\left( x\right) \right]
-1\right) \right. \\ && \left. \mbox{}  +\sum_{i=1}^{J}n\mu _{i}u_{i}\left( \exp \left[ -nV^{n}\left( x+
\frac{1}{n}\pi (x,\tilde{v}_{i})\right) +nV^{n}\left( x\right) \right]
-1\right) - nc\right] 
\end{eqnarray*}
results. \ Recall the definition $l(x)=x\log x-x+1$ for $x>0$. We now
divide throughout by $n$ and
use the convex duality relation 
\[
\left[ e^{y}-1\right] =\sup_{x>0}[xy-l(x)]
\]
to represent the terms in the previous display. \ For example, in the sum on 
$j$ we take $x=\bar{\lambda}_{j}/\lambda _{j}$ and $y=-\left[
nV^{n}\left( x+\frac{1}{n}v_{j}\right) -nV^{n}\left( x\right) \right] $.
Representing each term in this way and multiplying by $-1$ produces the
first line in (\ref{eq:dpe}). \ The boundary condition that is the second line in (\ref{eq:dpe})
follows directly from the relation between $W^{n}$ and $V^{n}$.
\qed

\noi{\bf Proof of Lemma \ref{lem:unif}:}
 We reduce the Lipschitz property on $(n^{-1}\Z^J_+)\cap\bar G$
to a Lipschitz property near the boundary.
To this end we use the following coupling.
  For $z\in G^n$, let $u^n(z)$ be a minimizer
in (\ref{eq:dpe3}). Given a point $x$ on the lattice, let $X^x$ denote
the process corresponding to the generator $\calL^{n,u^n}$ and starting at $x$
(see the discussion following (\ref{eq:gens})). 
To simplify the notation we will not explicitly denote the dependence of quantities such as $X^x$ on $n$.
Let $u(t)=u^n(X^x(t))$, and let $F_t$ be the filtration generated by $X^x$.

Fix a point $ y \not= x $ and let $ X^y $ denote the queueing process
on this probability space that 
starts at $ y $ and  uses the control $u$.
In other words,
$X^y$ is the image,
under the Skorokhod map,
of $y+X^x(\cdot)-x$.
The evolution of the processes $X^x$ and $X^y$ are identical, save that jumps which would cause $X^y$ to leave $\Z^J_+$ are deleted.
Automatically, $u$ is suboptimal for the control problem starting from $y$.
 Define
\be
  V^n(y ; \policy ) = - n^{-1}\log E_x^{\policy ,n} e^{-nc\sig^y }
\ee
where $ \sig^y $ is the exit time of $  X^y $ from $G$.
Note that due to the coupling we may take expectations with respect
to $ E_x^{u,n} $ rather then with respect to $ E_y^{u,n} $.
Since $(X^y,u)$ is a (possibly suboptimal) controlled Markov process, we have
\be\label{eq:VbyVu}
V^n (x ) - V^n ( y ) \leq V^n (x) - V^n (y ; \policy ) .
\ee

Define $ \sig = \min \{ \sig^x , \sig^y \} $.
By Theorem~\ref{th:SP} on the Lipschitz continuity
of the Skorokhod map we have
\be
\dist (  X^x (\sig ) , \partial G ) \leq K_1 | x - y | , \quad
\dist (  X^y (\sig ) , \partial G ) \leq K_1 | x - y |,
\ee
since at least one of the processes has left $G$ by $ \sig $.
In the last display, $K_1$ is the constant appearing in (\ref{sp:lip}). We claim that
\be\label{dist:1}
V^n (x ) - V^n ( y ) \leq \sup \{ V^n (z) : z \in S \} \ee
where $ S \doteq \{z\in n^{-1}\Z_+^J\cap \bar G:
 \dist (z , \pl_{co} G ) \leq K_1 |x-y| \} $.
To establish this, note that
\beaa
V^n (x ) - V^n ( y , \policy ) 
  &=& - \frac 1 n \left [ \log E_x^{u,n} e^{-nc \sigma^x } 
                              - \log E^{u,n}_x e^{-nc \sigma^y } \right ] \\
 & \leq & - \frac 1 n \left [ \log E_x^{u,n} \left [ e^{-nc \sigma }
        E_x^{u,n} \left ( e^{-nc (\sigma^x - \sigma )} \left |  X^x (\sigma)
           \right .  \right ) \right ] - \log E_x^{u,n} e^{-nc \sigma} \right ] \\
 & \leq & \sup_{z\in S} - \frac 1 n \left [ \log E_x^{u,n} \left [ e^{-nc \sigma }
                 E_x^{u,n} \left ( e^{-nc (\sigma^x - \sigma )} 
               \left |  X^x (\sigma) = z \right .  \right ) \right ]
              - \log E_x^{u,n} e^{-nc \sigma} \right ]\\
& =     & \sup_{z\in S} - \frac 1 n \left [ \log \left [ E_x^{u,n} 
             \left ( e^{-nc (\sigma^x - \sigma )} \left |  X^x (\sigma)  = z 
            \right .  \right ) E_x^{u,n} e^{-nc \sigma } \right ]
              - \log E_x^{u,n} e^{-nc \sigma} \right ]\\
 & = & \sup_{z\in S} - \frac 1 n \left [ \log E_x^{u,n} \left ( 
           e^{-nc (\sigma^x - \sigma )} \left |  X^x (\sigma)  = z 
               \right .  \right ) \right ]
\eeaa
However, by the strong Markov property,
$$
- \frac 1 n \left [ \log E_x^{u,n} \left ( e^{-nc (\sigma^x - \sigma )} 
                \left |  X^x (\sigma)  = z \right .  \right ) \right ]
\leq \sup_u - \frac 1 n \log \left ( E_z^{u,n} e^{-nc \sigma^z }\right ) 
= V^n (z)
$$
and together with~(\ref{eq:VbyVu}) we have~(\ref{dist:1}).

To prove the lemma, one needs to show that
$|V_n(x)-V_n(y)|\le c_0|x-y|$ for all $n$ and all
$x,y\in(n^{-1}\Z_+^J)\cap\bar G$, where $c_0$ does not depend
on $x,y$ and $n$. It suffices to prove this inequality for
$x,y$ such that $|x-y|=n^{-1}$.
Since the roles of $x $ and $y$ are symmetric, and in view of
(\ref{dist:1}), it suffices to show that
for $\{x\in G:\dist(x,\pl_{co} G)\le K_1n^{-1}\}$,
\begin{equation}\label{eq:55}
 V_n(x)=-n^{-1}\log\inf_uE_x^{u,n}e^{-nc\sig^x}\le c_1n^{-1},
\ee
where $c_1 > 0$ is a constant.

Let us first treat the case where $G$ is not a rectangle.
In that case, Condition \ref{cond:G} implies that for
any $x$ with $\dist(x,\pl_{co} G)\le K_1n^{-1}$,
\be\label{imp:1}
\mbox{ there is } i\in\calJ_+\ \mbox{ such that }\ x+c'n^{-1}e_i\not\in G,
\ee
where $c'$ is a constant.
Let such $i$ be fixed.
To show (\ref{eq:55}), it is enough to show that for any
$x$ such that $\dist(x,\pl_{co}G)\le K_1n^{-1}$, and any
$n$ and $u$,
\begin{equation}\label{eq:56}
 E_x^{u,n}e^{-nc\sig^x}\ge c_2>0.
\ee
Recall that $\la_i>0$. Let $S_t$ denote the event that
all service processes and all arrival processes, except for the
one corresponding to $i$, do not increase on $[0,t]$.
Recall that the expected time till a Poisson process of rate $\la$
hits level $K$ is $K/\la$.
Then for any $\al\in(0,1)$
\beaa
 E_x^{u,n}e^{-nc\sig^x} & \ge & \al P_x^{u,n}(e^{-nc\sig^x}>\al)\\
 &=&
\al P_x^{u,n}\lt(\sig^x<-\frac{\log\al}{nc}\rt).
\eeaa
Choosing $t_0=-(\log\al)/nc = 2\tilde c/n\la_i$
and using $P_x^{u,n}(\sig^x<2E\sig^x)\ge1/2$,
\beaa
 E_x^{u,n}e^{-nc\sig^x} & \ge & \al P_x^{u,n}(\sig^x<t_0|S_{t_0})P_x^{u,n}(S_{t_0})\\
 &\ge&
e^{-2 c \tilde c / \la_i } \frac 1 2 c_3
\eeaa
where $c_3 > 0 $ is the probability that a Poisson process with rate
$ n c_4 $ has not jumped by time $ t_0 = 2 \tilde c / \la_i n $.
This proves (\ref{eq:56}), which implies (\ref{eq:55}),
and hence the statement of the lemma holds.

In the case where $G$ is a rectangle,
the bound~(\ref{dist:1}) does not suffice since $V^n (x) $ is
discontinuous near $\pl_cG$.
We therefore prove that a similar bound applies, where there
supremum is over $S=\{z\in(n^{-1}\Z_+^J)\cap\bar G:
\dist(z,\pl_oG)\le K_1|x-y|\}$.
To apply the previous argument we need to show that if $  X^x (t)$ is
close to $\pl_cG$,
then neither $  X^x $ nor $ X^y $ will exit (locally) through that boundary.
This is clear for $  X^x $: the only
way for the process to leave $G$ is due to a service to one of the queues,
say, queue $j$, leading to an increase in queue $i$.
However, allowing this service is certainly not optimal: it is better to
avoid this control, as our objective is to increase $ \sigma^x $.
To prevent $  X^y $ from exiting we need to modify the coupling argument as
follows. The control $u^y$ used by $ X^y$ avoids a jump that leads
$  X^y $ to exit through $\pl_cG$ (that is, queue $j$ above will not be
served if $  X_i^y (t) = z_i - 1 $.) Note that this is the only possible
type of jump that leads the process out of $G$.
Moreover, the $ \ell_1 $ distance
between $  X^x $ and $ X^y $ may only decrease due to this change in
control: the control is changed only if $  X^x_i (t) < X^y_i (t) $,
and following the service $  X^x_i $ increases by $1$ so that
$ |  X^x_i (t) -  X^y_i (t) | $ decreases by $1$, while
$  X^x_j $ decreases by $1$.

Condition \ref{cond:G} still
implies (\ref{imp:1}), but only for $x$ such that
$$
\mbox{$z_i-x_i\le n^{-1}$, for some $i\in\calJ_+$.}
$$
For such $x$, the argument in the last paragraph holds.
However, for $x$ near $\pl_c G$
there is nothing to prove, since the process never exits through such a
boundary.
\qed

\noi{\bf Proof of Lemma \ref{lem:bound:m}:}
The first part is an immediate consequence of the fact that
 $ u_i \geq 0 $ and both
$\calL^{n,u,m}V^n(x)$ and $\rho(u,m)$ depend on $u$ as $\sum_iu_i\eta_i$,
where $\eta_i$ is a function of $m_i ,x,n$ but not of $u$.

For the second part of the lemma, one can explicitly solve for
$m^n$ in terms of $V^n$, and get $m^n(x,u)=m^n(x)=
((\bar\la^n_i(x)),(\bar\mu^n_i(x)))$,
where
$$
 \bar\la^n_i(x)=\la_ie^{-n\del_i V^n(x)}, \quad
 \bar\mu^n_i=\mu_ie^{-n\tilde\del_i V^n(x)},
$$
and
$$
 \del_i V^n(x)\doteq V^n(x+n^{-1}v_i)-V^n(x), \quad
 \tilde\del_i V^n(x)\doteq V^n(x+n^{-1}\pi(x,\tilde v_i))-V^n(x).
$$
The result follows from Lemma \ref{lem:unif},
since it shows that there is a constant $b_2$
independent of $x,n$ where
$$
 n\del_i V^n(x) \ge-b_2, \quad
 n\tilde\del_i V^n(x) \ge-b_2.
$$
\qed

\noi{\bf Proof of Lemma \ref{lem:cont}:}
We fix $b\in [b^*,\infty]$ and suppress it from the notation throughout
the proof.
Item 1 of the lemma is trivial under Condition \ref{cond:G}.1.
Under Condition \ref{cond:G}.2, by continuity of the functions
$ \phi_i$, we only need to show is that $\pl_{co}G_a$ and $\pl_{co}G$
do not intersect. Consider first $a>0$, and let $x\in\pl_{co}G_a$.
Then $x_i=a+ \phi_i(x_1,\ldots,x_{i-1},x_{i+1},\ldots,x_J)$
for some $i\in\calJ_+$, and therefore
$x$ cannot belong to the closure of $G$. The proof for $a<0$ is
similar.

Let $\beta_0[u](t)=m_0$ for all $u,t$, where $m_0$ sets all $\bar\la_i=
\la_i$ and $\bar\mu_i=0$. Then $\rho(u(t),\beta_0[u](t))$ is bounded
by a constant, and the dynamics, unaffected by $u$, follow $ X(t)=
x+\sum_i\la_ie_it$ and leave the bounded set $G$ within a finite time bounded
by ${\rm diam}(G)/\max_i\la_i<\infty$. Therefore
\[
V^- (x)\le \sup_u C (x,u,m_0) = C(x,\beta_0) \le c_1< + \infty .
\]
Similarly,
\[
V^+ (x)\le \sup_\alpha C (x,\alpha (m_0),m_0) \le c_1 < + \infty .
\]

It is useful to notice that for all $a\in(0,a_0)$
and $y\in\pl_o G$ there is $i=i_y\in\calJ_+$ such that
$y+2ae_{i_y}\not\in G_a$. Similarly, for all $a\in(-a_0,0)$
and $y\in\pl_o G_a$ there is $i=i_y\in\calJ_+$ such that
$y+2ae_{i_y}\not\in G$.

First consider $a>0$ and recall that $ \sigma_a $ (resp., $\sigma $)
is the exit time from $G_a $ (resp., $G$), so that for any fixed $ u$ and
$ \beta [u] $ we have $ \sigma \leq \sigma_a $.
Therefore, since $c$ and $ \rho $ are positive,
\beaa
 V^-_a(x)
 &=    & \inf_\beta\sup_u \int_0^{\sig_a} (c+\rho(u(s), \beta[u](s))ds \\
 &\geq & \inf_\beta\sup_u \int_0^{\sig} (c+\rho(u(s), \beta[u](s))ds \\
 & = & V^- (x) .
\eeaa
Thus to prove the Lipschitz property a one-sided bound suffices.
Recall that $ \phi (\sigma ) $ is the exit point from $ G$ and
for each $\beta$ define the extension $\beta_a$ by
$$
 \beta_a[u]=\lt\{\begin{array}{ll}\beta[u](t) & t\in[0,\sig), \\ 
 \hat m & t\in[\sig,\infty),
 \end{array}\rt.
$$
where $\hat m$ sets all $\bar\mu_j=0$ and
$\bar\la_j=1_{j=i_{\phi(\sig)}}$.
Then for any $\beta$
\beaa
 V^-_a(x)
 &=& \inf_{\beta}\sup_u C_a(x,\beta[u],u) \\
 &\le& \inf_{\beta}\sup_u C_a(x,\beta_a[u],u)\\
 &=& \inf_\beta\sup_u\lt[C(x,\beta[u],u)+\int_\sig^{\sig_a}
 (c+\rho(u(s),\hat m)ds\rt] \\
 &\le&
V^-(x)+c_1a,
\eeaa
where the last line follows since $\rho(u(s),\hat m)$ is bounded
and since by the previous paragraph, $\sig_a-\sig\le 2a$.
Note that $c_1$ does not depend on $b\in[b^*,\infty]$.
For $a<0$ the same argument shows that $V^-(x)\le V^-_a(x)+c_3|a|$,
by interchanging the roles of $G$ and $G_a$.

For $ V_a^+ (x) $ note that an argument as above gives 
$ V_a^+ (x) \geq V^+ (x) $. For each $m$ define $ m_a $ by
$$
 m_a = \lt\{ \begin{array}{ll} m (t) & t\in[0,\sig), \\ 
 \hat m & t\in[\sig,\infty),  
 \end{array}\rt .
$$
where $ \hat m $ is as above. Let $ \alpha_\epsilon $ be an $\epsilon
$-optimal strategy. Then, since for any fixed $ u$ and
$ m $ we have $ \sigma \leq \sigma_a $,
\beaa
 V^+_a(x)
 &=& \sup_\alpha \inf_m C_a (x,m, \alpha[m]) \\
 & \leq & \inf_m C_a (x,m, \alpha_\epsilon [m]) + \epsilon \\
 & \leq & \inf_m C_a (x,m_a , \alpha_\epsilon [m_a ]) + \epsilon,
\eeaa
since we are taking the infimum over a smaller class of controls.
By the definition of $ C_a $,
\beaa
V^+_a(x)
 & \leq & \inf_m C (x, m , \alpha_\epsilon [m ])
           + \int_\sigma^{\sigma_a}
               (c+\rho(\alpha_\epsilon [m_a ](s),\hat m (s))\, ds + \epsilon \\
 & \leq & \sup_\alpha \inf_m C (x,m, \alpha[m]) + c_2 a + \epsilon
\eeaa
by the previous argument,
where $c_2 $ does not depend on $x$, $ \epsilon $ and $b$.
Since $ \epsilon $ is
arbitrarily small, the proof for $a>0 $ and $ V_a^+ (x) $ is established.
\qed

\noi{\bf Proof of Lemma \ref{lem:apriori}:}
We suppress $b$ from the notation, throughout the proof.
It is obvious that one can restrict the infimum
over $\beta\in B$ to the class of strategies
$\beta$ for which $C(x,\beta)\le V^-(x)+1$.
Within this class, for every $\beta$ and $u$,
$$
\sigma c \leq \int_0^\sig[c+\rho(u(t),\beta[u](t))]dt\le V^-(x)+1,
$$
and therefore one always has that $\sig\le T_0\doteq(V(x)+1)/c$.
Lemma \ref{lem:apriori} asserts an upper bound on the cost till time a fixed time $T_0$,
and so we must define the strategy for times $t \in [\sigma, T_0]$.
Let $\hat m$ be an arbitrary fixed element of $M$.
Then the extended $\beta$ is just 
$$
 \hat\beta[u](t)=\lt\{\begin{array}{ll}\beta[u](t) & t<\sig.\\ 
 \hat m & t\ge\sig,
 \end{array}\rt.
$$
With this definition
one has that $C(x,u,\beta[u])=C(x,u,\hat\beta[u])$. One can
therefore further restrict to strategies $\beta$ satisfying
$\hat\beta=\beta$. For such $\beta$, it follows that
$$
 \int_0^{T_0}\rho(u(t),\beta[u](t))dt \le c_1T_0,
$$
where $c_1$ does not depend on $u,\beta$.
The result regarding $V^-$ follows.

Regarding $V^+$,
let $m_0$ be a control which sets all $\mu_i$ and $\la_i$ to
zero, except that $\la_{i_0}=1$ for some $i_0\in\calJ_+$.
Then for any $\al\in A$ and $m$ for which
\begin{equation}\label{eq:20}
C(x,\al[m],m)\le
C(x,\al[m_0],m_0),
\end{equation}
one has $c\sig(x,\al[m],m)
\le C(x,\al[m],m)\le
C(x,\al[m_0],m_0)\le c_1<\infty$. Note that $c_1$ can be chosen
independent of $\al$, since the dynamics and running cost
under $m_0$ are independent of $\al$. Clearly, for each $\al$ it suffices
to consider, in optimizing over $m$, only those $m$ that satisfy
(\ref{eq:20}). It follows that it suffices to consider only
those $m$ for which $\s(x,\al[m],m)\le c_1/c$.
This completes the proof of the lemma.
\qed

\noi{\bf Proof of Lemma \ref{lem:vallip}:}
Fix $ b\in[b^*,\infty]$ which we omit from the notation.
Recall from Lemma \ref{lem:cont} that $V^\pm$ are bounded
on $G$.
We first show that $V^-$ is Lipschitz. Assume first that
Condition \ref{cond:G}.2 holds.
Recall that for $x\in G$,
$$
 V^-(x)=\inf_\beta\sup_u C(x,\beta[u],u).
$$
Let $\beta_\eps^x$ be an $\eps$-optimal strategy starting
from $x$, i.e.,
$$
\sup_u C(x,\beta_\eps^x[u],u)\le V^-(x)+\eps.
$$
For any $z\in G$ let $\s_z=\inf\{t:\phi_z\not\in
G\}$, where $\phi_z$ is the solution to
$\dot\phi=\pi(\phi,v(u,\beta_\eps^x[u]))$, with $\phi(0)=z$.
Note that $C(x,u,\beta_\eps^x[u])=\int_0^{\s_x}[c+\rho(u(t),
\beta_\eps^x[u](t))]dt$ (with possibly $ \sigma_x = \infty $).
Now let $y\in G$. Note that on $[0,\s_x\w\s_y]$, one has by
the Lipschitz property of the Skorokhod map
that $|\phi_x (t) -\phi_y (t)|\le c_1|x-y|$,
where $c_1$ is some constant. 
Recall that we are considering the case of Condition \ref{cond:G}.2.
Therefore, at $\s_x\w \s_y$, both
$\phi_x$ and $\phi_y$ are within a distance of $c_1|x-y|$ of
the boundary $\pl_o G$. Because of the assumptions on the domain $G$,
there exists a constant $c_2$
 such that at time $\s_x\w\s_y$, both
$\phi_x+c_2e_{i^*}\not\in G$ and
$\phi_y+c_2e_{i^*}\not\in G$, where $i^*\in\calJ_+$
and, moreover, $c_2 $ is
independent of $ x,y, b \geq b^* $ and $ i^* $.

Define $\beta_\eps^{x,y}$ as
$\beta_\eps^{x,y}[u]=\beta_\eps^x[u]$ on $[0,\sig_x)$, and,
if $\sig_y\ge\sig_x$, set $\beta_\eps^{x,y}[u]=m_0$ on $[\s_x,\s_y]$.
Here, $m_0$ sets all $\bar\mu_i$ and all $\bar\la_i$ to zero, except
that it sets $\bar\la_{i^*}=\la_{i^*}$, where $i^*$ is as above.
Consequently, $ \sig_y < \sig_x + c^\prime | x - y | $ for some 
$ c^\prime > 0 $, and there exists $u_\eps$ such that
\beaa
 V^-(y) &\le& \sup_u C(y,u,\beta_\eps^{x,y}[u]) \\
 &\le&
\int_0^{\s_y}[c+\rho(u_\eps(s),\beta_\eps^{x,y}[u_\eps](s))]ds+\eps\\
 &\le&
\int_0^{\s_x}[c+\rho(u_\eps(s),\beta_\eps^{x,y}[u_\eps](s))]ds
 +1_{\s_x<\s_y}
 \int_{\s_x}^{\s_y}[c+\rho(u_\eps(s),\beta_\eps^{x,y}[u_\eps](s))]ds
 +\eps\\
 &\le&
C(x,u_\eps,\beta_\eps^{x}[u_\eps])+c_3|x-y|+\eps\\
 &\le&
\sup_u C(x,u,\beta_\eps^x[u])+c_3|x-y|+\eps\\
 &\le&
V^-(x)+c_3|x-y|+2\eps.
\eeaa
Since $x,y\in G$ and $\eps>0$ are arbitrary, and $c_3$ does not
depend on them or on $ b$, $V^-$ is Lipschitz, uniformly for
$ b\in[b^*,\infty] $.

In case that Condition \ref{cond:G}.1 holds, the same
argument shows that $V^-_a(y)\le V^-(x)+c_3|x-y|+2\eps$,
where $a=c_1|x-y|$. By Lemma \ref{lem:cont}, this implies that
$V^-(y)\le V^-(x)+c_4|x-y|+2\eps$, some constant $c_4$, and therefore
$V^-$ is Lipschitz.

Next, consider the upper value
$$
 V^+(x)=\sup_\al\inf_mC(y,\al[m],m)
$$
under Condition \ref{cond:G}.2.
Let $x,y\in G$. Note that there is an $\al_\eps^x$ such that
$$
 V^+(x)\le\inf_m C(x,\al_\eps^x[m],m)+\eps.
$$
and an $m_\eps=m_\eps(x,y)$ for which
\beaa
V^+(y) &\ge& \inf_m C(y,\al_\eps^x[m],m) \\
 &\ge& C(y,\al_\eps^x[m_\eps],m_\eps)-\eps.
\eeaa
Let $\s_z=\inf\{t:\phi_z\not\in G\}$, where $\phi_z$ is the solution
to $\dot\phi=\pi(\phi,v(\al_\eps^x[m_\eps],m_\eps))$, with
$\phi(0)=z$. Let $i^*$ be defined in an analogous way to that in
the first paragraph of the proof. Now
define $\bar m_\eps=\bar m_\eps(x,y)$ as follows. If
$\sig_x\le\sig_y$, let $\bar m_\eps=m_\eps$. If $\sig_y<\sig_x$, let
$\bar m_\eps$ agree with $m_\eps$ on $[0,\sig_y)$ and with
$m_0$ on $[\sig_y,\s_x]$. Here, $m_0$ sets all $\bar\mu_i$
and all $\bar\la_i$ to zero, except that it sets $\bar\la_{i^*}
=\la_{i^*}$.
Since $m_\eps$ and $\bar m_\eps$ agree on $[0,\s_y)$, the restrictions
to $[0,\sig_y]$ of $\al_\eps^x[m_\eps]$ and of $\al_\eps^x[\bar m_\eps]$
agree a.e.\ on $[0,\sig_y]$, and therefore,
$$
 C(y,\al_\eps^x[m_\eps],m_\eps)=C(y,
 \al_\eps^x[\bar m_\eps],\bar m_\eps).
$$
Arguing again by the Lipschitz property of the
Skorokhod map and the definition of $m_0 $,
there is a constant $c_4$ for which
$(\sig_x-\sig_y)^+\le c_4|x-y|$. Hence
\beaa
V^+(y) &\ge& C(y,\al_\eps^x[\bar m_\eps],\bar m_\eps)-\eps \\
&\ge & \int_0^{\sig_x}[c+\rho(\al_\eps^x[\bar m_\eps^x](s),\bar m_\eps^x(s))]
ds - 1_{\sig_y<\sig_x}\int_{\sig_y}^{\sig_x}
[c+\rho(\al_\eps^x[\bar m_\eps^x](s),\bar m_\eps^x(s))]ds -\eps\\
&\ge& C(x,\al_\eps^x[\bar m_\eps],\bar m_\eps)-c_5|x-y|-\eps\\
&\ge& \inf_m C(x,\al_\eps^x[m],m)-c_5|x-y|-\eps\\
&\ge& V^+(x)-c_5|x-y|-2\eps.
\eeaa
Since $c_5$ does not depend on $x,y,\eps$ or $b$, we have that $V^+$
is Lipschitz uniformly for $b\in[b^*,\infty] $.

Under Condition \ref{cond:G}.1, the same argument shows that
$V^+_a(y)\ge V^+(x)-c_5|x-y|-2\eps$, where $a=c_1|x-y|$, and
again one argues by Lemma \ref{lem:cont}.
\qed

\noindent{\bf Proof of Lemma \ref{lem:xbar}:}
The processes are constructed recursively using a sequence
of standard exponential clocks.
Recall that $\calL^{n,u,m}$ is given for every $n$, $u\in U$,
$m\in M$ by
\[
\calL^{n,u,m}f(x) =
\sum_{j=1}^J n \bar \lambda_j [f(x+{n}^{-1} v_j)-f(x)]
+ \sum_{i=1}^J n \bar \mu_i u_i [f(x+n^{-1}\pi(x,\tilde v_i))-f(x)].
\]
Given $n$, $x_n$ and $\beta$, we construct a filtered
probability space and three processes, $\bar X(t)$, $\bar u(t)$
and $\bar m(t)$ (to simplify notation, we do not write the superscript
$n$ in the notation of $\bar X^n$, $\bar u^n$ and $\bar m^n$)
such that (a) $\bar X,\bar u$ and $\bar m$
are $(\bar F_t)$-adapted;
(b) $\bar m(t)=\beta[\bar u](t)$ a.e.\ $t\ge0$, a.s.;
(c) $\bar u(\cdot )=u^n(X(\cdot ))$ a.s.\ (where $u^n$ is as in
the statement before the lemma); and (d) for any $f$, the process
$$
f(\bar X(t))-\int_0^t\calL^{n,\bar u(s),\bar m(s)}f(\bar X(s))ds
$$
is an $(\bar F_t)$-martingale.
For (a--d) to hold, it suffices that (a--c) hold, and (e)
on any finite interval the process $\bar X$ jumps finitely many 
times---we denote
the $k$th jump by $\tau_k$ and let $\tau_0=0$; (f) the random times
$(\tau_k)$ are
stopping times on $(\bar F_t)$, and (g) denoting $X_k=\bar X(\tau_k)$,
for any $k$,
$$
\bar E[f(X_{k+1})-f(X_k)|\bar F_{\tau_k}] =
 \sum_{i=1}^J \bar E[A_i^{k,\bar u,\bar m}
 +B_i^{k,\bar u,\bar m}|\bar F_{\tau_k}],
$$
where
$$
A_i^{k,\bar u,\bar m}= n\int_{\tau_k}^{\tau_{k+1}}
   \bar\la_i(s)ds [f(X_k+n^{-1}v_i)-f(X_k)],
$$
$$
B_i^{k,\bar u,\bar m} = n\int_{\tau_k}^{\tau_{k+1}} \bar\mu_i(s)u_i(s)ds
    [f(X_k+n^{-1}\pi(X_k,\tilde v_i))-f(X_k)].
$$

The construction is recursive.  On a complete probability
space $(\bar\Om,\bar F,\bar P)$ we are given $2J$ independent
i.i.d.\ standard Poisson processes, denoted
$a_i$ and $b_i$, $i=1,\ldots,J$. Let $T^a_i(k)$ [resp., $T^b_i(k)$]
denote the first time $a_i$ [resp., $b_i$] equals $k$.
For each $\om\in\Om$ we construct recursively a sequence of times
$(\tau_k)$ and the processes $\bar X, \bar u$ and $\bar m$
up to time $\tau_k$. Once these processes are
defined, we will define $(\bar F_t)$, $\bar F_t\subset \bar F$,
$t\ge0$, and verify
that items (a--c), (e--g) are satisfied on
$(\bar\Om,\bar F,(\bar F_t),\bar P)$.

We set $\bar X(0)=x_n$ and $\bar u(0)=u^n(x_n)$.
Since $\bar m$ need only be defined
almost everywhere on $[0,\infty)$, we do not define it at zero nor
at any $\tau_k$, $k=1,2,\ldots$.
Now assume that we have constructed $\tau_i$, $i\le k$ as well as
the processes $\bar X$ and $\bar u$ on
$[0, \tau_k]$ and $\bar m$ a.e.\ on $[0,\tau_k]$.
Let
$\hat u^k(t)=u^n(\bar X(t\w \tau_k))$, $t\ge0$. Let also
$\hat m^k=\beta[\hat u^k]$.
With $\hat u^k(\cdot )=(\hat u_i^k(\cdot ))$ and
$\hat m^k(\cdot )=((\hat\la_i^k(\cdot )),(\hat\mu_i^k(\cdot )))$, let
$$
p_i^k(t)=n\int_0^t\hat\la_i^k(s)ds,
\qquad q_i^k(t)=n\int_0^t\hat\mu_i^k(s)\hat u_i^k(s)ds,
\qquad i=1,\ldots,J, \ t\ge0.
$$
Denoting $\Del z(s)=z(s)-z(s-)$, let also
$$
 \tau_{k+1}=\inf\{t>\tau_k: \mbox{ either $\Del a_i(p_i^k(t))>0$
 or $\Del b_i(q_i^k(t))>0$ for some $i=1,\ldots,J$}\},
$$
where $\inf\emptyset=+\infty$.
We first consider the case that $\tau_{k+1}<+\infty$.
In this case,
\be\label{eq:exi}
\mbox{there is $i$ such that either $\Del a_i(p_i^k(\tau_{k+1}))>0$ or
$\Del b_i(q_i^k(\tau_{k+1}))>0$.}
\ee
In the former case we let
$\hat v^k= v_i$; otherwise we let $\hat v^k=\tilde v_i$.

The three processes are defined on the next interval
as follows. Let $\bar X(t)=\bar X(\tau_k)$ for  $t\in(\tau_k,\tau_{k+1})$,
and $\bar X(\tau_{k+1})=\bar X(\tau_k)+n^{-1}\pi(\bar X(\tau_k),\hat v^k)$.
Let $\bar u(t)=u^n(\bar X(t))$ for $t\in(\tau_k,\tau_{k+1}]$.
Let $\check u(t)=u^n(\bar X(t\w\tau_{k+1}))$ and define
$\bar m(t)=\beta[\check u](t)$, $t\in[0,\tau_{k+1}]$.
Note that since $\beta$ is a strategy,
this definition of $\bar m$ is consistent with its definition
up to $\tau_k$ since so is the definition of $\bar u$.
For the same reason, for a.e.\ $t\le\tau_{k+1}$,
$\bar m(t)=\hat m^k(t)$.
In particular, the equations for $p_i^k, q_i^k$ still hold
if we replace hats by bars, namely,
\be\label{eq:pq}
p_i^k(t)=n\int_0^t\bar\la_i(s)ds,
\quad q_i^k(t)=n\int_0^t\bar\mu_i(s)\bar u_i(s)ds,
\quad i=1,\ldots,J, \ \tau_k\le t\le \tau_{k+1}.
\ee
Note that the above relations are consistent in the sense that
for a given $k$, they hold not only for
$t\in[\tau_k,\tau_{k+1}]$, but in fact for $t\in[0,\tau_{k+1}]$.
Hence, on the event $\tau_k\to\infty$, one can equivalently
consider the processes
\be\label{def:pq}
p_i(t)=n\int_0^t\bar\la_i(s)ds,
\quad q_i(t)=n\int_0^t\bar\mu_i(s)\bar u_i(s)ds,
\quad i=1,\ldots,J, \ t\ge0.
\ee
This completes the definition of the three processes on $[0,\tau_{k+1}]$.

In case that $\tau_{k+1}=+\infty$, the definitions above of
$\bar X$, $\bar u$ and $\bar m$
all apply on $(\tau_k,\tau_{k+1})$ and there is
nothing else to define.

To complete the construction of the three processes on
$\bar\Om\times[0,+\infty)$, we must consider the set $\Om_0$
of $\om\in\bar\Om$
for which $\bar T\doteq\sup\tau_k$ is finite.
We show that this set is $\bar P$-null
owing to the fact that the range $\bar M^b $ of $\beta$ consists
of bounded functions.
Suppose $\bar T$ is finite. The construction above defines $\bar X,\bar u$
and $\bar m$
on $[0,\bar T)$. Let $\bar u'(t)=\bar u(t)$ for $t<\bar T$
and define $\bar u'(t)$ arbitrarily
on $[\bar T,+\infty)$ but such that $\bar u'\in\bar U$. Then $\bar m'
=\beta[\bar u']$
agrees with $\bar m$ a.e.\ on $[0,\bar T]$. Since each component of $\bar m'$
is bounded by $b$,
\be\label{eq:boundsum}
 n^{-1}\max_{i=1}^J [p_i(\bar T)\vee q_i(\bar T)]\le
 2J\bar T b <+\infty.
\ee
However, by construction, $\bar T<\infty$ implies that
either $a_i(p_i(t))\to\infty$ or
$b_i(q_i(t))\to\infty$ as $t\uparrow\bar T$, for some $i$.
Hence $\bar T<\infty$ must be a null set.
We let $\bar X,\bar u$ and $\bar m$ be defined arbitrarily on $\Om_0$.

The definition of the process $\bar Y$ is similar to that of
$\bar X$, but where $\pi(x,v)$ is replaced by $v$ throughout.
The relation $\bar X=\Gamma(\bar Y)$ is clear from the construction.

Define for each $t\ge0$ $\bar F_t$ to be the $\sig$-field
generated by $\{ \bar Y(s), s\in[0,t] \} $. Note that it is
equivalently defined as the $\sig$-field generated by
$\{a_i(p_i(t)), b_i(q_i(t)), i=1,\ldots,J\}$, where
$p_i,q_i$ are as in (\ref{def:pq}).
By construction, $\bar u(t)=u^n(\bar X(t))$,
$t\ge0$ and item (c) holds. Item (b), namely that $\bar m=\beta[\bar u]$,
also holds by construction.
$\bar X$ and $\bar u$ are therefore $(\bar F_t)$-adapted,
and since $\beta$ is a strategy, so is $\bar m$, and item (a) holds.
Items (e) and (f) are trivial. Concerning (g),
let $i^k\in\{1,\ldots,2J\}$ denote the index $i$ satisfying
(\ref{eq:exi}) in case that $\Del a_i(p_i^k(t))>0$ holds,
and let it denote $i+J$ in the case $\Del b_i^k(q_i^k(t))>0$.
It suffices to show that for every $i\in\{1,\ldots,2J\}$,
$$
 \bar P(i^k=i|\bar F_{\tau_k})=
 \lt\{\begin{array}{ll}
 \bar E[\int_{\tau_k}^{\tau_{k+1}}p_i(s)ds|\bar F_{\tau_k}]/Z_k
  & i\le J,\\[.1in]
 \bar E[\int_{\tau_k}^{\tau_{k+1}}q_{i-J}(s)ds|
 \bar F_{\tau_k}]/Z_k & i> J,
 \end{array}\rt.
$$
where $Z_k$ is a normalization factor (not depending on $i$).
For $k=0$ ($\tau_k=0$),
this is a well known property of exponential clocks.
For $k>0$, the same argument holds, merely because conditional on
$\bar F_{\tau_k}$, the processes $\int_{\tau_k}^\cdot p_i(s)ds$,
$\int_{\tau_k}^\cdot q_i(s)ds$ are independent, and moreover,
$a_i(\cdot-\tau_k)-a_i(\tau_k),b_i(\cdot-\tau_k)-b_i(\tau_k)$
are still independent Poisson processes (which is a statement
on the lack of memory for exponential random variables).

The proof of the claim regarding the martingale associated with
$\calL_0$ is similar (only simpler). This completes the proof of
the first part of the lemma.

Clearly,
$$
   \max_i p_i(T_0)\vee q_i(T_0)\le nT_0b,
$$
where $T_0$ is as in Lemma \ref{lem:apriori}.
Thus, if $N_n=\max\{k:\tau_k\le T_0\}$, then
$$
 N_n\le \sum_i a_i(nT_0b)+b_i(nT_0b),
$$
and (\ref{ineq:Nn}) follows.
\qed

\noi{\bf Proof of Lemma \ref{lem:xbar2}:}
The proof is completely
analogous to that of Lemma \ref{lem:xbar}, and is therefore
omitted.
\qed

\noi{\bf Proof of Theorem \ref{th:b}:}
By Theorems \ref{th:unique} and \ref{th:solve},
$V^{b,-}=V^{b,+}$ for all $b\in[b^*,\infty)$. As a result,
Theorem \ref{th:limit}, implies that
$V^n\to V^{b,-}$ for all $b\in[b^*,\infty)$, as $n\to \infty$.
In particular,
$V^{b,-}$ does not depend on $b\in[b^*,\infty)$.
It remains to show that for all $x$, $V^{b,-}(x)\to V^-(x)$
and $V^{b,+}(x)\to V^+(x)$ as $b\to\infty$.

\noi\uu{\em Proof that $V^{b,-}\to V^-$.}
It is immediate from the definitions that $V^-\le V^{b,-}$.

Let $\beta\in B$, and let $\s=\s(x,u,\beta)$ be the exit
time of $\phi$ from $G$ where $\dot\phi=\pi(\phi,v(u,\beta[u]))$,
$\phi(0)=x$. Let $\bar\beta$ be defined by
$$
\bar\beta[u](t)=\lt\{\begin{array}{ll}\min\{b,\beta[u](t)\}
 & t\le\s, \\  \hat m & t>\s,
\end{array}\rt.
$$
where $\hat m$ sets all $\bar\mu_j=0$ and $\bar\la_j=1_{j=i_{\phi(\s)}}$,
and the minimum is componentwise. It is clear that
$\bar\beta$ is a strategy. Let $\bar u$ be any extension
of $u$ to $[0,\infty)$, and denote by $\bar\phi$ and $\bar\s$ the
dynamics and exit time corresponding to $x,\bar\beta,\bar u$. 
Recall that by (\ref{eq:lal}) $b$ is greater than all $\lambda_i$ and $\mu_i$.
Thus
\beaa
 C(x,\bar u,\bar \beta[\bar u]) &=&
\int_0^\s(c+\rho(u,\bar\beta[u])ds
 +1_{\bar\s>\s}\int_\s^{\bar \s}(c+\rho(\bar u,\hat m))ds\\
 &\le&
C(x,u,\beta[u])+c_1(\bar\s-\s)^+.
\eeaa
Moreover, by the Lipschitz property of the Skorokhod map, and denoting
$\beta[u]=(\bar\la_i,\bar\mu_i)$,
$$
(\bar\s-\s)^+
\le c_2|\phi(\s)-\bar\phi(\s)|\le
c_2\int_0^\s\sum_i[(\bar\la_i-b)^+ +(u_i\bar\mu_i-b)^+]ds
$$
Since it is enough to consider $\beta$ for which (for any $u$)
$\la_i$ and $u_i\mu_i$ are uniformly integrable over $[0,\s]$,
we have that $(\bar\s-\s)^+\le \del(b)$, where $\del(b)\to0$
as $b\to\infty$. This shows that
$\lim_{b\to\infty}V^{b,-}(x)\le V^-(x)$.

\noi\uu{\em Proof that $V^{b,+}\to V^+$.}
It is immediate that $V^+\le V^{b,+}$.

To show that $V^+(x)\ge \lim_{b\to\infty}V^{b,+}(x)$
it is enough to show that for $b\ge b^*$ and $a$ small,
$V^+(x)\ge V^{b,+}_{-a}(x)$.
For any $m\in\bar M$ let $m^b$ denote the pointwise and componentwise
truncation of $m$ at level $b$. For any $\al\in A$, let $\al^b\in A$
be defined by $\al^b[m]=\al[m^b]$.
We will write $m\in\bar M(\al,a)$ if $m,\al,a$ satisfy
$C_a(x,\al[m],m)\le V^+_a(x)+1$.
In the expression for $V^+_a(x)$,
$$
 \sup_\al\inf_m C_a(x,\al[m],m),
$$
it is enough to consider $\al\in A$ and $m\in\bar M(\al,a)$
(including for $a=0$).
For such $\al,m$, the functions
$\bar\la_i,u_i\bar\mu_i$ are uniformly integrable over $[0,T]$.
Let $\al\in A$ and $m\in\bar M(\al^b,0)$. Consider
a truncation of $\bar\la_i$ and $\bar\mu_i$ at $b$.
Denote by $\phi$ [resp., $\phi^b$] the dynamics that correspond to
$(x,\al^b,m)$, [resp., $(x,\al^b,m^b)$].
Then the effect of the truncation
on $\phi$ is such that for all $a>0$ there is $b$ such
that $\sup_{[0,T]}|\phi-\phi^b|\le a$ (by uniform integrability).
In particular, $|\phi^b(\s)-\phi(\s)|\le a$.
Hence, using
the monotonicity of the running cost for large values of
the rates, and that $\al^b[m^b]=\al^b[m]$,
$$
 C_{-a}(x,\al^b[m^b],m^b)\le C(x,\al^b[m],m).
$$
We thus have
$$
C_{-a}(x,\al[m^b],m^b)\le C(x,\al^b[m],m).
$$
Since $m\in\bar M(\al^b,0)$ implies that $m^b\in\bar M(\al^b,-a)$,
\beaa
 \inf_{m\in\bar M^b}C_{-a}(x,\al[m],m) &\le&
\inf_{m:m^b\in\bar M(\al^b,-a)}
 C_{-a}(x,\al[m^b],m^b)
 \\ &\le& \inf_{m\in\bar M(\al^b,0)}C(x,\al^b[m],m)
 \\ &=&
\inf_{m\in\bar M}C(x,\al^b[m],m).
\eeaa
Hence
$$
\sup_{\al\in A}\inf_{m\in\bar M^b}C_{-a}(x,\al[m],m)\le
\sup_{\al\in A}\inf_{m\in\bar M}C(x,\al[m],m).
$$
Taking $a\to0$ by letting $b\to\infty$, we have from
Lemma \ref{lem:cont} that $\lim_bV^{b,+}(x)\le V^+(x)$.
\qed

\vspace{\baselineskip}
\noindent
{\it 1991 Mathematics Subject Classification.}
Primary 60F10, 60K25;
Secondary 93E20, 60F17.

\bibliographystyle{plain}

\end{document}